# 2nd order PDEs: geometric and functional considerations

## *Second order partial derivatives equations that can be solved by the characteristics method*

Andrea Pezzi

195756@edu.unito.it , pezzi.andrea@gmail.com

## *Abstract*

*This papers deals with partial differential equations of second order, linear, with constant and not constant coefficients, in two variables, which admit real characteristics. I face the study of PDEs with the mentality of the applied physicist, but with a weakness for formalization: look inside the* black box *of the formulas, try to compact them (for example, proceeding from an inverse transformation of coordinates) and make them* smart *(in the context, reformulating the theory by means of differential operators and related invariants), applying them with awareness and then connecting them to geometry or to spatial categories, which are in mathematics what is closest to the sensible reality. Finally, proposing examples that are exercise and corroborating for theory.*

## *Topics*

*The geometric meaning of "invariant to a differential operator".*

*Operator "Principal Part" and its factorization: commutativity and product with and without "residues"(first order terms). Related conditions by operators and invariants derivatives.*

*Coordinate transformation by invariants and expression of the hyperbolic and parabolic operators in the new coordinates.*

*Properties of the Jacobian Matrix and relations between invariants derivatives and inverse coordinates transformation or the initial variables derivatives.*

*Commutativity conditions and product without "residues" in terms of inverse coordinate transformations that allow to build commutative differential operators or whose product is without residues (or both).*

*Diffeomorphisms and plane transformations: new operators and invariants in the new coordinate space which lead to the "chain rule" in compact form.*

*Conclusive considerations and examples who compares different methods of solution.*

# PDE II: considerazioni geometriche e funzionali

## *Equazioni alle derivate parziali del secondo ordine risolvibili con il metodo delle caratteristiche*

Andrea Pezzi

195756@edu.unito.it , pezzi.andrea@gmail.com

### ***Introduzione***

La presente trattazione riguarda le equazioni differenziali alle derivate parziali del secondo ordine, lineari, a coefficienti costanti e non costanti, in due variabili, che ammettono caratteristiche reali. Affronto lo studio delle PDE con la mentalità del fisico applicato che però ha un debole per la formalizzazione: guardare cosa succede dentro la scatola nera delle formule, cercare di compattarle (procedendo ad esempio da una trasformazione inversa di coordinate) e renderle più *smart* (nel contesto riformulando la teoria mediante gli operatori differenziali e i relativi invarianti), applicandole però con consapevolezza e quindi collegandole alla geometria ovvero alle categorie spaziali, che sono in matematica ciò che più si avvicina alla realtà sensibile. Infine, proponendo esempi che siano di esercizio e corroborante per la teoria.

### ***Argomenti trattati***

Il significato geometrico di *invariante ad un operatore differenziale*.

Operatore *Parte Principale* e la sua fattorizzazione: commutatività e prodotto con e senza *residui* (termini del primo ordine). Condizioni relative attraverso operatori e derivate degli invarianti.

Trasformazione di coordinate per mezzo degli invarianti ed espressione degli operatori iperbolico e parabolico nelle nuove coordinate.

Proprietà della Matrice Jacobiana e relazione tra le derivate degli invarianti e la trasformazione di coordinate inversa, ovvero le derivate delle variabili iniziali.

Condizioni di commutatività e prodotto senza residui in termini della trasformazione inversa di coordinate, le quali consentono di costruire operatori differenziali commutativi oppure il cui prodotto sia senza residui (o entrambe le cose).

Diffeomorfismi e trasformazioni del piano: nuovi operatori e invarianti nelle nuove coordinate che conducono alla *chain rule* in forma compatta.

Considerazioni conclusive ed esempi che confrontano diversi metodi di soluzione

## *Sulle PDE del secondo ordine con caratteristiche reali*

Sia

$$au_{xx} + 2bu_{xy} + cu_{yy} + du_x + eu_y + gu = f \quad (*)$$

un'equazione differenziale di secondo grado alle derivate parziali, con

$$u = u(x,y), \ \Delta = b^2 - ac > 0 \quad \textbf{(problema iperbolico)}$$

Per il momento si consideri la *parte principale* $au_{xx} + 2bu_{xy} + cu_{yy}$ e l'operatore differenziale $\mathcal{L}$

$$\mathcal{L}[u] = au_{xx} + 2bu_{xy} + cu_{yy}$$

Ricordando la scomposizione del trinomio di secondo grado, se i coefficienti *a, b, c* sono costanti si può senz'altro scrivere:

$$\mathcal{L} = a\left(\frac{\partial}{\partial x} - \Lambda^+ \frac{\partial}{\partial y}\right)\left(\frac{\partial}{\partial x} - \Lambda^- \frac{\partial}{\partial y}\right) = a\,\mathcal{L}^+\mathcal{L}^-$$

con $\Lambda^\pm = \frac{-b \pm \sqrt{\Delta}}{a}$ .

Se si ha un operatore differenziale nella forma

$$\mathcal{L} = \left(\alpha \frac{\partial}{\partial x} + \beta \frac{\partial}{\partial y}\right) \quad con \quad \begin{cases} \alpha = \alpha(x,y) \\ \beta = \beta(x,y) \end{cases} \neq 0$$

e si vuole trovare una funzione $\phi(x,y): \mathcal{L}[\phi] = 0$ , ovvero *invariante* rispetto all'applicazione dell'operatore $\mathcal{L}$:

$$\mathcal{L}[\phi] = \alpha\phi_x + \beta\phi_y = 0 \ \Rightarrow \ \frac{\phi_x}{\phi_y} = -\frac{\beta}{\alpha} \quad (**)$$

Dalla proporzione (\*\*), differenziando si ha che:

$$d\phi = \phi_x\,dx + \phi_y\,dy = K\,(\beta\,dx - \alpha\,dy)$$

dove *K* è un fattore di proporzionalità, la cui natura ci si propone ora di indagare.

Si cominci col prendere *K=1*, $d\phi = \beta\,dx - \alpha\,dy$.

Si ponga ora $\tilde{\phi} = F(\phi)$; differenziando segue :

$$d\tilde{\phi} = F'(\phi)\,d\phi = F'(\phi)\,(\beta\,dx - \alpha\,dy)$$

per cui $K = F'(\phi)$ e si conclude che un invariante non è univocamente determinato: se $\phi$ è tale per cui $\mathcal{L}[\phi] = 0$, allora qualsiasi funzione che abbia $\phi$ come argomento sarà altrettanto invariante rispetto a $\mathcal{L}$.

*Es. Si consideri l'operatore:*

$$\mathcal{D} = (x+y)\frac{\partial}{\partial x} + (x-y)\frac{\partial}{\partial y}$$

*Integrando la forma differenziale*

$$d\phi(x,y) = (x-y)dx - (x+y)dy$$

*si ottiene subito che*

$$\phi(x,y) = \frac{x^2}{2} - xy - \frac{y^2}{2}$$

*è invariante rispetto all'applicazione dell'operatore $\mathcal{D}$. Ma non è l'unico: lo è qualsiasi funzione di $\phi$.*

Sia ora $\phi(x,y) = k$, $k \in \mathbb{R}$, una funzione implicita formata dalla famiglia di curve di livello della superficie $z = \phi(x,y)$. Allora:

$$d\phi = \phi_x\, dx + \phi_y\, dy = dk = 0$$

Posto $\phi_y \neq 0$ ($\phi$ è un fascio di curve regolari),

$$\frac{\phi_x}{\phi_y} = -\frac{dy}{dx} \Rightarrow \frac{dy}{dx} = \frac{\beta}{\alpha}$$

Si tratta di fatto di risolvere una PDE I omogenea.

*Es. Sia*

$$\mathcal{R} = -y\frac{\partial}{\partial x} + x\frac{\partial}{\partial y}$$

$$\frac{dy}{dx} = \frac{\beta}{\alpha} = -\frac{x}{y} \Rightarrow ydy = -xdx \Rightarrow \phi(x,y) = x^2 + y^2 = k$$

*Ponendo $k = r^2$ si hanno le circonferenze di centro l'origine e raggio r; considerato che $\partial_x$, $\partial_y$ sono versori dello spazio tangente, l'operatore $\mathcal{R}$ induce una rotazione per la quale $\phi(x,y)$, formata dalla collezione di tali curve è evidentemente invariante.*

Tornando all'operatore $\mathcal{L}^\pm$, siccome $a(\Lambda^\pm)^2 + 2b(\Lambda^\pm) + c = 0$, $-\Lambda^\pm$ verificherà l'uguaglianza $a(\Lambda^\pm)^2 - 2b(\Lambda^\pm) + c = 0$ e poiché $\frac{dy}{dx} = -\Lambda^\pm$, si può scrivere l'*equazione alle caratteristiche*:

$$a\left(\frac{dy}{dx}\right)^2 - 2b\left(\frac{dy}{dx}\right) + c = 0 \quad (***)$$

dalla quale si ricavano gli invarianti $\phi(x,y) = k$, $\psi(x,y) = h$, ad $\mathcal{L}^+$ ed $\mathcal{L}^-$ rispettivamente.

E' utile alla soluzione delle PDE II fare alcune considerazioni di carattere geometrico.

$\phi(x,y) = k$ è una funzione implicita che è anche la maniera più generale di descrivere un fascio di curve cartesiane $y = f(x)$ al variare di $k$.

E' abbastanza intuitivo il fatto che il vettore gradiente $\nabla\phi$ sia normale alle linee di livello della superficie $z = \phi(x,y)$, rappresentate proprio da $\phi(x,y) = k$. Inoltre $\mathcal{L} = \left(\alpha\frac{\partial}{\partial x} + \beta\frac{\partial}{\partial y}\right)$ può interpretarsi come un prodotto scalare $\mathcal{L} = \boldsymbol{v} \cdot \boldsymbol{\nabla}$, con $\boldsymbol{v} = (\alpha, \beta)$, o ancora come *derivata direzionale "non normalizzata"* $\frac{\partial}{\partial v} = \boldsymbol{v} \cdot \boldsymbol{\nabla} = \|v\|\frac{\partial}{\partial \hat{v}}$, $\hat{v} = \frac{v}{\|v\|}$ (nel contesto ciò che è rilevante è la direzione più che il modulo). Se $\phi$ è invariante rispetto all'operatore $\mathcal{L}$ allora la sua derivata lungo la direzione $\boldsymbol{v}$ è nulla. Visto che $\nabla\phi = \boldsymbol{n}$, con $\boldsymbol{n}$ vettore normale alle curve $\phi(x,y) = k$, poiché per ipotesi $\nabla\phi \cdot \boldsymbol{v} = \boldsymbol{n} \cdot \boldsymbol{v} = \boldsymbol{0}$, segue che $\boldsymbol{v}$ è il vettore tangente alle curve $\phi(x,y) = k$. In sostanza, cercare un invariante all'operatore differenziale $\mathcal{L} = (\alpha\partial_x + \beta\partial_y)$, consiste nel trovare una superficie $z = \phi(x,y)$ le cui curve di livello (la loro proiezione) coincidano con le linee di campo del vettore $\boldsymbol{v} = (\alpha, \beta)$, in modo tale che la derivata direzionale $\mathcal{L}[\phi] = \frac{\partial\phi}{\partial v}$ sia nulla, essendo per definizione z costante lungo le curve $\phi(x,y) = k$. Tali curve, in questo contesto, possono essere chiamate *curve caratteristiche*.

*Es. Si riconsideri* $\mathcal{R} = -y\,\partial_x + x\,\partial_y$. *Il vettore $\boldsymbol{v}$ ha componenti $(-y, x)$ e si è già calcolato* $\phi(x,y) = x^2 + y^2$, *infatti* $\mathcal{R}[\phi] = -y(2x) + x(2y) = 0$. *Si confrontino i grafici:*

*I* 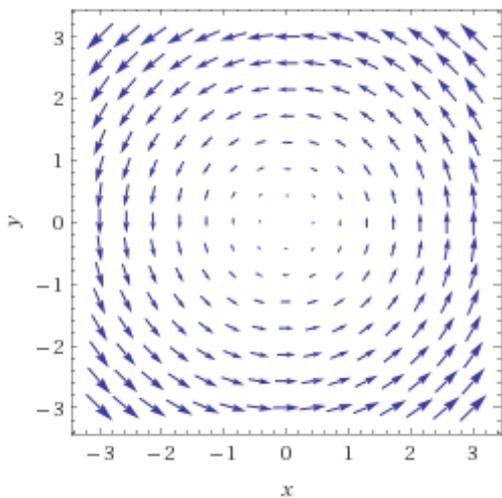 *II* 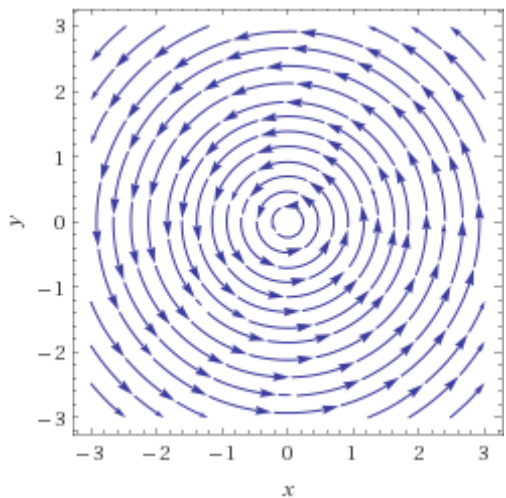

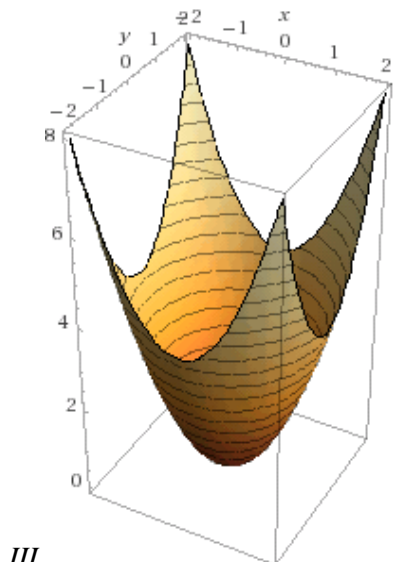

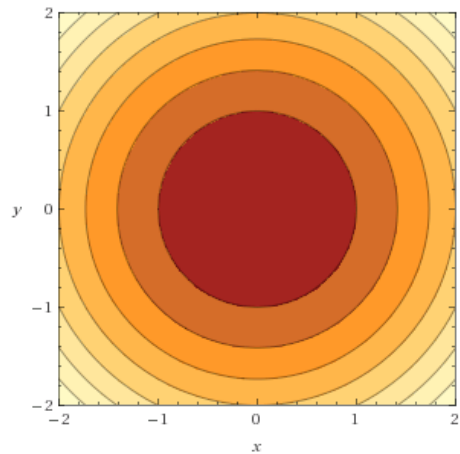

*III* *IV*

*Il grafico I è il campo di direzioni $\boldsymbol{v}=(-y,x)$ lungo cui si deriva; nel II sono rappresentate le linee di campo del vettore $\boldsymbol{v}$, ovvero quelle linee che hanno proprio $\boldsymbol{v}$ come vettore tangente in ogni punto. Il III è la superficie $\phi(x,y) = x^2 + y^2$ con le curve di livello $\phi(x,y) = x^2 + y^2 = k$ nel IV. Si può notare come II e IV coincidano.*

In questo caso, ricordando la (**):

$$\mathcal{L}^\pm = \boldsymbol{w}^\pm \cdot \boldsymbol{\nabla} = \frac{\partial}{\partial \boldsymbol{w}^\pm} \quad , \quad \boldsymbol{w}^\pm = (1, -\Lambda^\pm)$$

$$\mathcal{L}^+[\phi] = (\phi_x, \phi_y) \cdot \left(1, -\frac{\phi_x}{\phi_y}\right) = 0$$

$$\mathcal{L}^-[\psi] = (\psi_x, \psi_y) \cdot \left(1, -\frac{\psi_x}{\psi_y}\right) = 0$$

$\partial u / \partial \boldsymbol{w}^\pm$ si può dunque interpretare come derivata di $u(x,y)$ nella direzione tangente a $\phi(x,y) = k$ ($\psi(x,y) = h$). Il fatto che $\mathcal{L}^-\mathcal{L}^+[u]$ misuri la variazione di $u(x,y)$ prima lungo le curve $\phi = k$ e poi lungo le curve $\psi = h$ (o viceversa), suggerisce che un cambiamento di coordinate in cui si usino proprio tali linee possa agevolare l'integrazione di (*).

Ricapitolando, una volta risolta l'equazione (***) e integrate le due ODE, $\frac{dy}{dx} = -\Lambda^\pm$, si ottengono due famiglie di caratteristiche $\phi(x,y) = \xi, \psi(x,y) = \eta$ che saranno le nuove coordinate, lasciando intese $\phi$ e $\psi$ come mappe di transizione.

Sia dunque $u = u(x,y)$ una funzione di classe $\mathcal{C}^2$ definita su un aperto $\Omega \subseteq \mathbb{R}^2$ e $\phi(x,y), \psi(x,y)$ diffeomorfismi in $\Omega$, per cui $x = \phi^{-1}(\xi, \eta), y = \psi^{-1}(\xi, \eta)$. Siano $\mathcal{L}^\pm$ due operatori differenziali nelle coordinate $(x,y)$ tali per cui $\mathcal{L}^+[\phi] = 0$, $\mathcal{L}^-[\psi] = 0$.

Si consideri la funzione $u(x,y)$ nelle nuove coordinate, $U(\xi, \eta) = U(\phi(x,y), \psi(x,y))$. Si ha che:

$$\mathcal{L}^-[U] = \frac{\partial U}{\partial \xi} \mathcal{L}^-[\phi] + \frac{\partial U}{\partial \eta} \mathcal{L}^-[\psi] = \frac{\partial U}{\partial \xi} \mathcal{L}^-[\phi]$$

$$\mathcal{L}^+[U] = \frac{\partial U}{\partial \xi} \mathcal{L}^+[\phi] + \frac{\partial U}{\partial \eta} \mathcal{L}^+[\psi] = \frac{\partial U}{\partial \eta} \mathcal{L}^+[\psi]$$

per cui si ha la corrispondenza:

$$\frac{\partial}{\partial \xi} = \frac{\mathcal{L}^-}{\mathcal{L}^-[\phi]} \qquad \frac{\partial}{\partial \eta} = \frac{\mathcal{L}^+}{\mathcal{L}^+[\psi]}$$

Si definisce commutatore l'operatore:

$$[\mathcal{L}^-, \mathcal{L}^+] = \mathcal{L}^- \mathcal{L}^+ - \mathcal{L}^+ \mathcal{L}^-$$

Se i due operatori commutano, cioè se $\mathcal{L}^- \mathcal{L}^+ = \mathcal{L}^+ \mathcal{L}^-$, segue che $[\mathcal{L}^-, \mathcal{L}^+] = 0$.

Ora, anche se in generale i due operatori non commutano, $[\mathcal{L}^-, \mathcal{L}^+] \neq 0$, l'uguaglianza delle derivate miste $\partial_\xi \partial_\eta = \partial_\eta \partial_\xi$ implica che

$$[\frac{\mathcal{L}^-}{\mathcal{L}^-[\phi]}, \frac{\mathcal{L}^+}{\mathcal{L}^+[\psi]}] = 0$$

Sempre procedendo dalla derivate miste si calcola:

$$\frac{\partial}{\partial \xi}\frac{\partial}{\partial \eta} = \frac{\mathcal{L}^-}{\mathcal{L}^-[\phi]}\left[\frac{\mathcal{L}^+}{\mathcal{L}^+[\psi]}\right] = \frac{\mathcal{L}^-\mathcal{L}^+}{\mathcal{L}^-[\phi]\,\mathcal{L}^+[\psi]} - \frac{\mathcal{L}^-(\mathcal{L}^+[\psi])}{\mathcal{L}^-[\phi]\,\mathcal{L}^+[\psi]}\frac{\mathcal{L}^+}{\mathcal{L}^+[\psi]}$$

$$\frac{\partial^2}{\partial \xi \partial \eta} = \frac{1}{\mathcal{L}^-[\phi]\,\mathcal{L}^+[\psi]}\left\{\mathcal{L}^-\mathcal{L}^+ - \mathcal{L}^-(\mathcal{L}^+[\psi])\frac{\partial}{\partial \eta}\right\}$$

da cui si trova la corrispondenza

$$\mathcal{L}^-\mathcal{L}^+ = \mathcal{L}^-[\phi]\,\mathcal{L}^+[\psi]\frac{\partial^2}{\partial \xi \partial \eta} + \mathcal{L}^-(\mathcal{L}^+[\psi])\frac{\partial}{\partial \eta}$$

Ripetendo lo stesso conto per $\partial_\eta \partial_\xi$ si ricava

$$\mathcal{L}^+\mathcal{L}^- = \mathcal{L}^-[\phi]\,\mathcal{L}^+[\psi]\frac{\partial^2}{\partial \xi \partial \eta} + \mathcal{L}^+(\mathcal{L}^-[\phi])\frac{\partial}{\partial \xi}$$

Nel caso più semplice in cui i due operatori commutino, è possibile scambiare l'ordine di applicazione degli operatori nei secondi addendi che risulteranno nulli:

$$\mathcal{L}^-(\mathcal{L}^+[\psi]) = \mathcal{L}^+(\mathcal{L}^-[\psi]) = 0 \qquad \mathcal{L}^+(\mathcal{L}^-[\phi]) = \mathcal{L}^-(\mathcal{L}^+[\phi]) = 0$$

e alla fine resta:

$$\mathcal{L}^-\mathcal{L}^+ = \mathcal{L}^-[\phi]\,\mathcal{L}^+[\psi]\frac{\partial^2}{\partial \xi \partial \eta} \qquad (°)$$

Nella (°) è presente solo una derivata mista, mentre le altre derivate seconde si sono annullate.

Visto che $\mathcal{L}^\pm = \left(\frac{\partial}{\partial x} - \Lambda^\pm \frac{\partial}{\partial y}\right)$ e $\Lambda^+ = \frac{\phi_x}{\phi_y}$, $\Lambda^- = \frac{\psi_x}{\psi_y}$, si ha

$$\mathcal{L}^-\mathcal{L}^+[U] = -(\phi_y \psi_y)^{-1}\left|\frac{\partial(\phi,\psi)}{\partial(x,y)}\right|^2 \frac{\partial^2 U}{\partial \phi \partial \psi} = -\phi_y \psi_y \begin{vmatrix} 1 & \Lambda^+ \\ 1 & \Lambda^- \end{vmatrix}^2 \frac{\partial^2 U}{\partial \phi \partial \psi}$$

Se i coefficienti di (*) non sono costanti ma dipendono da $x$ e $y$, si devono operare delle sostituzioni mediante le trasformazioni inverse, $x = \phi^{-1}(\xi,\eta)$ e $y = \psi^{-1}(\xi,\eta)$.

Infine, essendo $\Lambda^\pm = \frac{-b \pm \sqrt{\Delta}}{a}$ e nel caso in cui $\mathcal{L} = a\,\mathcal{L}^+\mathcal{L}^-$, si conclude che:

$$\mathcal{L}[U] = -\frac{4\Delta}{a}\,\phi_y \psi_y\,\frac{\partial^2 U}{\partial \xi \partial \eta}$$

Sommando i termini di grado inferiore si avrà (*) nella forma:

$$\frac{4\Delta}{a} \phi_y \psi_y \frac{\partial^2 U}{\partial \xi \partial \eta} = \mathcal{F}(\xi, \eta, U_\xi, U_\eta, U)$$

Per cui la *forma canonica* della (*) è del tipo

$$U_{\xi\eta} = F(\xi, \eta, U_\xi, U_\eta, U)$$

Come ci si aspettava i due operatori $\mathcal{L}^\pm$ nelle coordinate usuali sono stati sostituiti da due derivate parziali nelle nuove coordinate: $\mathcal{L}^- \mathcal{L}^+[u(x,y)] \xrightarrow{(\phi,\psi)} \partial_\xi \partial_\eta U(\xi,\eta)$.

*Es. Caso Omogeneo*

$$U_{\xi\eta} = 0$$

*L'integrale generale sarà*

$$U(\xi,\eta) = F(\xi) + G(\eta)$$

*Infatti derivando:* $U_\xi = \partial_\xi [F(\xi) + G(\eta)] = f(\xi)$ , $U_{\xi\eta} = \partial_\eta f(\xi) = 0$

*Se i coefficienti di (*) sono tali per cui (*) è omogenea e $\Lambda^\pm$ sono numeri reali:*

$$u(x,y) = F(y + \Lambda^+ x) + G(y + \Lambda^- x)$$

Bisogna ricordarsi che la (°) è valida nell'ipotesi $[\mathcal{L}^+, \mathcal{L}^-] = 0$, altrimenti potrebbe risultare incompleta o comunque dipendente dall'ordine di derivazione. E' possibile stabilire a priori che tale condizione sia verificata solo nel caso in cui i coefficienti di (*) siano costanti per cui si calcola, a partire dagli operatori:

$$\mathcal{L} = a\partial_{xx} + 2b\partial_{xy} + c\partial_{yy} \quad \text{ed} \quad \mathcal{L}^\pm = (\partial_x - \Lambda^\pm \partial_y)$$

$$\mathcal{L}^- \mathcal{L}^+ = a^{-1} \mathcal{L} - \mathcal{L}^-[\Lambda^+] \partial_y$$

$$\mathcal{L}^+ \mathcal{L}^- = a^{-1} \mathcal{L} - \mathcal{L}^+[\Lambda^-] \partial_y$$

e il commutatore da origine ad un nuovo operatore:

$$[\mathcal{L}^-, \mathcal{L}^+] = (\mathcal{L}^+[\Lambda^-] - \mathcal{L}^-[\Lambda^+]) \partial_y$$

che si annichilisce se:

$$[\mathcal{L}^-, \mathcal{L}^+] = 0 \iff \mathcal{L}^+[\Lambda^-] = \mathcal{L}^-[\Lambda^+]$$

che è la condizione di commutatività per cui è vera (°).

Dovrà inoltre essere:

$$\mathcal{L}^+[\Lambda^-] = \mathcal{L}^-[\Lambda^+] = 0 \iff \mathcal{L} = a \mathcal{L}^- \mathcal{L}^+ = a \mathcal{L}^+ \mathcal{L}^-$$

*Es.1*

$$u_{xy} + 2xu_{yy} = u_y$$

$$\mathcal{L} = \frac{\partial}{\partial y}\left(\frac{\partial}{\partial x} + 2x\frac{\partial}{\partial y}\right)$$

*In questo caso si possono tranquillamente scambiare gli operatori $\mathcal{L}^+ = \partial_x + 2x\partial_y$ ed $\mathcal{L}^- = \partial_y$, difatti:*

$$\mathcal{L}^\pm = (\partial_x - \Lambda^\pm \partial_y) \Rightarrow \Lambda^+ = -2x, \Lambda^- = -1 \Rightarrow \mathcal{L}^+[\Lambda^-] = \mathcal{L}^-[\Lambda^+] = 0 \Rightarrow [\mathcal{L}^-, \mathcal{L}^+] = 0$$

*Gli invarianti rispetto a $\mathcal{L}^\pm$ sono $\phi(x,y) = y - x^2 = \xi$ e $\psi(x) = x = \eta$ e si verifica facilmente:*

$$\mathcal{L}^+[\phi] = -2x + 2x = 0 \qquad \mathcal{L}^-[\psi] = \partial_y x = 0$$

Si calcoli ancora $\mathcal{L}^+[\psi] = \partial_x x + 2x\partial_y x = 1$ ; $\mathcal{L}^-[\phi] = \partial_y(y - x^2) = 1$ *e dalla (°) si ha:*

$$\mathcal{L}[u] = \mathcal{L}^- \mathcal{L}^+[u] = \frac{\partial^2 U}{\partial \xi \partial \eta}$$

*Data la corrispondenza* $\{\mathcal{L}^- \mathcal{L}^+\}_{(x,y)} = \{\partial_\xi \partial_\eta\}_{(\xi,\eta)}$ *con* $\mathcal{L}^- = \partial_y$, *si deduce: $u_y = U_\xi$. Allora la forma canonica dell'equazione sarà:*

$$U_{\xi\eta} = U_\xi$$

*da cui la soluzione:*

$$(U_\xi)_\eta = U_\xi \Rightarrow \frac{dU_\xi}{U_\xi} = d\eta \Rightarrow U_\xi = f(\xi)e^\eta \Rightarrow U(\xi,\eta) = e^\eta F(\xi) + G(\eta)$$

*Infine si sostituiscono a $\xi$ e $\eta$, $\phi(x,y)$ e $\psi(x)$:*

$$u(x,y) = e^x F(y - x^2) + G(x)$$

*Es.2*

*Allo stesso risultato si può arrivare scrivendo:*

$$u_{xy} + 2xu_{yy} - u_y = \partial_y(u_x + 2xu_y - u) = 0 \Rightarrow u_x + 2xu_y - u = f(x)$$

*e cioè trasformando una PDE II in una del primo ordine: $u_x + 2xu_y = u + f(x)$. Si scrive il sistema:*

$$\frac{dx}{1} = \frac{dy}{2x} = \frac{du}{u + f(x)}$$

*E' immediato il calcolo del primo invariante: $\varphi_1 = y - x^2$, per cui le caratteristiche (la loro proiezione) saranno delle parabole.*
*Più problematica è la ricerca del secondo invariante, soluzione di: $(u + f(x))dx - du = 0$, che richiede l'intervento di un fattore integrante $\mu(x) = e^{-x}$ per chiudere la forma differenziale e infine trovare :*
$\varphi_2 = ue^{-x} + g(x)$.
*La soluzione è una combinazione funzionale degli integrali primi:*

$$ue^{-x} + g(x) = F(y - x^2) \Rightarrow u(x,y) = e^x F(y - x^2) + G(x)$$

*Es.3*

*L'equazione*

$$8\,sh^2\left(\frac{x}{2}\right) u_{xx} - 4\,e^{\frac{y}{2}}\,sh(x)\,u_{xy} + 2\,e^y\,u_{yy} + \left(2\,e^{\frac{y}{2}} + e^y\right) u_y = 0$$

*che a prima vista sembra molto complicata, ammette in realtà una soluzione relativamente semplice.*

*Risolvendo rispetto alla parte principale si trova:*

$$\Lambda^- = \frac{1}{2}\,e^{\frac{y-x}{2}}\,csch\left(\frac{x}{2}\right) = \frac{e^{\frac{y}{2}}}{e^x - 1}$$

$$\Lambda^+ = \frac{1}{2}\,e^{\frac{y+x}{2}}\,csch\left(\frac{x}{2}\right) = \frac{e^{x+\frac{y}{2}}}{e^x - 1}$$

*da cui i due operatori $\mathcal{L}^-$ e $\mathcal{L}^+$. Calcolando esplicitamente $\mathcal{L}^+[\Lambda^-]$ e $\mathcal{L}^-[\Lambda^+]$ si ha:*

$$\mathcal{L}^+[\Lambda^-] = \mathcal{L}^-[\Lambda^+] = -\frac{1}{8}\,e^{\frac{y}{2}}\,csch^2\left(\frac{x}{2}\right)\left(2 + e^{\frac{y}{2}}\right) \neq 0$$

*I due operatori commutano e di per sé questo costituisce una prima semplificazione.*

*Inoltre se si divide per a ( il coefficiente di $u_{xx}$ ) il coefficiente di $u_y$ si ottiene proprio $-\mathcal{L}^+[\Lambda^-]$ ($-\mathcal{L}^-[\Lambda^+]$), perciò l'ultimo termine sarà $-a\,\mathcal{L}^-[\Lambda^+]\,u_y$.*

*Poiché la parte principale è $\mathcal{L}[u]$, l'equazione può scriversi :*

$$\left(\mathcal{L} - a\,\mathcal{L}^-[\Lambda^+]\,\partial_y\right) u(x,y) = 0$$

*Ma l'operatore tra parentesi non è altro che a $\mathcal{L}^-\mathcal{L}^+$ da cui per la (°) la forma canonica $U_{\xi\eta} = 0$ con soluzione $U(\xi,\eta) = F(\xi) + G(\eta)$.*

*Come ultimo passo si cercano gli invarianti sulle caratteristiche per il cambio di coordinate: da $\frac{dy}{dx} = -\Lambda^\pm$ seguono in questo caso due ODE a variabili separabili da cui la mappa*

$$\phi(x,y) = \ln(e^x - 1) - 2e^{-\frac{y}{2}} = \xi \qquad \phi^{-1}(\xi,\eta) = \xi - \eta = x$$

$$\psi(x,y) = \ln(e^x - 1) - 2e^{-\frac{y}{2}} - x = \eta \qquad \psi^{-1}(\xi,\eta) = -2\ln\left\{\frac{1}{2}\left[\ln(e^{\xi-\eta} - 1) - \xi\right]\right\} = y$$

*e infine l'integrale generale*

$$u(x,y) = F\left(\ln(e^x - 1) - 2e^{-\frac{y}{2}}\right) + G\left(\ln(e^x - 1) - 2e^{-\frac{y}{2}} - x\right)$$

*Non era necessario calcolare la mappa inversa ai fini della risoluzione dell'equazione, tuttavia si osservi come è più "semplice" l'espressione di $\phi^{-1}$ rispetto a $\psi^{-1}$: non è un caso, è infatti una conseguenza del commutare degli operatori, come si chiarirà in seguito.*

Volendo occuparsi dei termini di grado inferiore a quello massimo, oltre che a trovare dei procedimenti che non tengano conto della natura dell'operatore della parte principale, ovvero se iperbolico, parabolico o ellittico (quest'ultimo caso non viene trattato in quanto apre il campo ai numeri complessi e alle funzioni armoniche, le quali meriterebbero un proprio approfondimento), né della natura della sua scomposizione, si enuncia la seguente regola:

date due carte $X(x_1, \ldots, x_n), Y(y_1, \ldots, y_n)$, dove $x_i$ e $y_i$ sono le componenti o coordinate, vale la *regola della catena*:

$$\frac{\partial}{\partial x_i} = \sum_{j=1}^{n} \frac{\partial y_j}{\partial x_i} \frac{\partial}{\partial y_j}$$

Allora date le carte $(x, y)$ e $(\xi, \eta)$ e la funzione di transizione scritta in componenti, $\Phi(\phi(x,y), \psi(x,y))$, si possono scrivere in forma compatta:

$$\frac{\partial}{\partial x} = \Phi^{-1} \circ \left( \frac{\partial \Phi}{\partial x} \cdot \boldsymbol{\nabla}_{\xi,\eta} \right) \qquad \frac{\partial}{\partial y} = \Phi^{-1} \circ \left( \frac{\partial \Phi}{\partial y} \cdot \boldsymbol{\nabla}_{\xi,\eta} \right)$$

con:

$$\frac{\partial \Phi}{\partial x} = (\phi_x, \psi_x), \quad \frac{\partial \Phi}{\partial y} = (\phi_y, \psi_y), \quad \boldsymbol{\nabla}_{\xi,\eta} = (\partial_\xi, \partial_\eta)$$

Applicando gli operatori a $u(x,y)$ e $U(\xi, \eta)$ si ottiene:

$$u_x = \phi_x U_\xi + \psi_x U_\eta \qquad\qquad u_y = \phi_y U_\xi + \psi_y U_\eta$$

Per le derivate del secondo ordine si hanno le uguaglianze:

$$\frac{\partial^2}{\partial x^2} = \Phi^{-1} \circ \left[ \left( \frac{\partial \Phi}{\partial x} \cdot \boldsymbol{\nabla}_{\xi,\eta} \right)^2 + \frac{\partial^2 \Phi}{\partial x^2} \cdot \boldsymbol{\nabla}_{\xi,\eta} \right]$$

$$\frac{\partial^2}{\partial y^2} = \Phi^{-1} \circ \left[ \left( \frac{\partial \Phi}{\partial y} \cdot \boldsymbol{\nabla}_{\xi,\eta} \right)^2 + \frac{\partial^2 \Phi}{\partial y^2} \cdot \boldsymbol{\nabla}_{\xi,\eta} \right]$$

$$\frac{\partial^2}{\partial x \partial y} = \Phi^{-1} \circ \left[ \left( \frac{\partial \Phi}{\partial x} \cdot \boldsymbol{\nabla}_{\xi,\eta} \right)\left( \frac{\partial \Phi}{\partial y} \cdot \boldsymbol{\nabla}_{\xi,\eta} \right) + \frac{\partial^2 \Phi}{\partial x \partial y} \cdot \boldsymbol{\nabla}_{\xi,\eta} \right]$$

La mappa inversa $\Phi^{-1} = (\phi^{-1}(\xi, \eta), \psi^{-1}(\xi, \eta))$ viene giustamente applicata dopo il calcolo esplicito del termine tra parentesi, altrimenti si potrebbe pensare che $\boldsymbol{\nabla}_{\xi,\eta}$ agisca sulle derivate di $\Phi$, mentre si tratta di operazioni su carte diverse.

L'applicazione di queste formule non risparmia ridondanti cancellazioni, per cui è bene arrivare ad esprimere l'operatore $\mathcal{L}$ nelle variabili *x, y* in una combinazione di derivate nelle coordinate $\xi, \eta$.

Si riprenda dunque l'espressione

$$\mathcal{L}^- \mathcal{L}^+ = a^{-1} \mathcal{L} - \mathcal{L}^-[\Lambda^+] \partial_y$$

ovvero

$$a^{-1}\mathcal{L} = \mathcal{L}^-\mathcal{L}^+ + \mathcal{L}^-[\Lambda^+]\partial_y$$

ma

$$\mathcal{L}^-\mathcal{L}^+ = \mathcal{L}^-[\phi]\,\mathcal{L}^+[\psi]\,\partial_{\xi\eta} + \mathcal{L}^-(\mathcal{L}^+[\psi])\partial_\eta$$

e ora si sa anche che

$$\partial_y = \phi_y\,\partial_\xi + \psi_y\partial_\eta$$

Mettendo tutto insieme si arriva a:

$$\mathcal{L} = \Phi^{-1}\circ\left\{a\left\{\mathcal{L}^-[\phi]\,\mathcal{L}^+[\psi]\,\partial_{\xi\eta} + [\mathcal{L}^-(\mathcal{L}^+[\psi]) + \mathcal{L}^-[\Lambda^+]\psi_y]\partial_\eta + (\mathcal{L}^-[\Lambda^+]\phi_y)\partial_\xi\right\}\right\}$$

dove $\Phi^{-1}$ sta ad indicare che i termini nelle variabili $x$ e $y$ vanno convertiti nelle nuove variabili $\xi,\eta$.

Ripetendo lo stesso procedimento per

$$\mathcal{L}^+\mathcal{L}^- = a^{-1}\mathcal{L} - \mathcal{L}^+[\Lambda^-]\partial_y$$

ovvero

$$a^{-1}\mathcal{L} = \mathcal{L}^+\mathcal{L}^- + \mathcal{L}^+[\Lambda^-]\partial_y$$

tenuto conto che

$$\mathcal{L}^+\mathcal{L}^- = \mathcal{L}^-[\phi]\,\mathcal{L}^+[\psi]\,\partial_{\xi\eta} + \mathcal{L}^+(\mathcal{L}^-[\phi])\,\partial_\xi$$

si arriva ad una espressione equivalente alla precedente:

$$\mathcal{L} = \Phi^{-1}\circ\left\{a\left\{\mathcal{L}^-[\phi]\,\mathcal{L}^+[\psi]\,\partial_{\xi\eta} + [\mathcal{L}^+(\mathcal{L}^-[\phi]) + \mathcal{L}^+[\Lambda^-]\phi_y]\partial_\xi + (\mathcal{L}^+[\Lambda^-]\psi_y)\partial_\eta\right\}\right\}$$

Confrontandole si nota che:

$$\mathcal{L}^-(\mathcal{L}^+[\psi]) + \mathcal{L}^-[\Lambda^+]\psi_y = \mathcal{L}^+[\Lambda^-]\psi_y$$

$$\mathcal{L}^+(\mathcal{L}^-[\phi]) + \mathcal{L}^+[\Lambda^-]\phi_y = \mathcal{L}^-[\Lambda^+]\phi_y$$

da cui si ha l'espressione finale per *l'operatore iperbolico della parte principale nelle coordinate* $(\xi,\eta)$

$$\boxed{\mathcal{L} = \Phi^{-1}\circ\left\{a\left[(\mathcal{L}^-[\phi]\,\mathcal{L}^+[\psi])\frac{\partial^2}{\partial\xi\partial\eta} + (\mathcal{L}^-[\Lambda^+]\,\phi_y)\frac{\partial}{\partial\xi} + (\mathcal{L}^+[\Lambda^-]\,\psi_y)\frac{\partial}{\partial\eta}\right]\right\}}$$

Inoltre, poiché

$$\mathcal{L}^-(\mathcal{L}^+[\phi]) = a^{-1}\mathcal{L}[\phi] - \mathcal{L}^-[\Lambda^+]\,\phi_y = 0 \quad\Rightarrow\quad \mathcal{L}[\phi] = a\,\mathcal{L}^-[\Lambda^+]\,\phi_y$$

$$\mathcal{L}^+(\mathcal{L}^-[\psi]) = a^{-1}\mathcal{L}[\psi] - \mathcal{L}^+[\Lambda^-]\psi_y = 0 \implies \mathcal{L}[\psi] = a\mathcal{L}^+[\Lambda^-]\psi_y$$

si può alternativamente scrivere

$$\mathcal{L} = \Phi^{-1} \circ \left[ a\left(\mathcal{L}^-[\phi]\,\mathcal{L}^+[\psi]\right)\frac{\partial^2}{\partial \xi \partial \eta} + \mathcal{L}[\phi]\,\frac{\partial}{\partial \xi} + \mathcal{L}[\psi]\,\frac{\partial}{\partial \eta} \right]$$

La validità dell'espressione può essere provata in diversi modi. Calcolando ad esempio il commutatore

$$[\mathcal{L}^-,\mathcal{L}^+] = (\mathcal{L}^+[\Lambda^-] - \mathcal{L}^-[\Lambda^+])\,\partial_y = (\mathcal{L}^+[\Lambda^-] - \mathcal{L}^-[\Lambda^+])\,\phi_y\,\partial_\xi + (\mathcal{L}^+[\Lambda^-] - \mathcal{L}^-[\Lambda^+])\,\psi_y\partial_\eta$$

dalle precedenti si ha che

$$\mathcal{L}^-(\mathcal{L}^+[\psi]) = (\mathcal{L}^+[\Lambda^-] - \mathcal{L}^-[\Lambda^+])\psi_y$$

$$-\mathcal{L}^+(\mathcal{L}^-[\phi]) = (\mathcal{L}^+[\Lambda^-] - \mathcal{L}^-[\Lambda^+])\phi_y$$

perciò

$$[\mathcal{L}^-,\mathcal{L}^+] = \mathcal{L}^-(\mathcal{L}^+[\psi])\partial_\eta - \mathcal{L}^+(\mathcal{L}^-[\phi])\,\partial_\xi$$

che coincide proprio con:

$$[\mathcal{L}^-,\mathcal{L}^+] = \mathcal{L}^-[\phi]\,\mathcal{L}^+[\psi]\,\partial_{\xi\eta} + \mathcal{L}^-(\mathcal{L}^+[\psi])\,\partial_\eta - \mathcal{L}^-[\phi]\,\mathcal{L}^+[\psi]\,\partial_{\xi\eta} - \mathcal{L}^+(\mathcal{L}^-[\phi])\,\partial_\xi$$

Inoltre, tenendo presente che $\Lambda^+ = \frac{\phi_x}{\phi_y}$, $\Lambda^- = \frac{\psi_x}{\psi_y}$, si può procedere ad un calcolo esplicito, prendendo ad esempio:

$$\mathcal{L}^-(\mathcal{L}^+[\psi]) = \frac{1}{\phi_y{}^2\psi_y}\{\phi_y{}^2(\psi_{xx}\psi_y - \psi_{xy}\psi_x) + \phi_x\,\phi_y(\psi_{yy}\psi_x - \psi_{xy}\psi_y)$$
$$- \psi_y{}^2(\phi_{xx}\,\phi_y - \phi_{xy}\,\phi_x) - \psi_x\psi_y(\phi_{yy}\phi_x - \phi_{xy}\phi_y)\}$$

$$\mathcal{L}^-[\Lambda^+]\psi_y = \frac{1}{\phi_y{}^2\psi_y}\{\psi_y{}^2(\phi_{xx}\phi_y - \phi_{xy}\phi_x) + \psi_x\psi_y(\phi_{yy}\phi_x - \phi_{xy}\phi_y)\}$$

$$\mathcal{L}^-(\mathcal{L}^+[\psi]) + \mathcal{L}^-[\Lambda^+]\psi_y =$$
$$= \frac{1}{\phi_y{}^2\psi_y}\{\phi_y{}^2(\psi_{xx}\psi_y - \psi_{xy}\psi_x) + \phi_x\,\phi_y(\psi_{yy}\psi_x - \psi_{xy}\psi_y)\} = \mathcal{L}^+[\Lambda^-]\psi_y$$

Queste espressioni potrebbero costituire il legame tra gli invarianti e le condizioni di commutatività, infatti $\mathcal{L}^-[\Lambda^+] = \mathcal{L}^+[\Lambda^-]$ si traduce in

$$\phi_y{}^2(\psi_{xx}\psi_y - \psi_{xy}\psi_x) + \phi_x\phi_y(\psi_{yy}\psi_x - \psi_{xy}\psi_y)$$
$$= \psi_y{}^2(\phi_{xx}\phi_y - \phi_{xy}\phi_x) + \psi_x\psi_y(\phi_{yy}\phi_x - \phi_{xy}\phi_y)$$

per cui si avrebbe giustamente $\mathcal{L}^-(\mathcal{L}^+[\psi]) = 0$. Dividendo per $\phi_y{}^2\psi_y{}^2$ e tenendo presente che $\frac{\phi_x}{\phi_y} = \Lambda^+$ e $\frac{\psi_x}{\psi_y} = \Lambda^-$

$$\frac{1}{\psi_y}\{\psi_{xx} - (\Lambda^+ + \Lambda^-)\psi_{xy} + \Lambda^+\Lambda^-\psi_{yy}\} = \frac{1}{\phi_y}\{\phi_{xx} - (\Lambda^+ + \Lambda^-)\phi_{xy} + \Lambda^+\Lambda^-\phi_{yy}\}$$

Ricordando che $\Lambda^\pm = \frac{-b \pm \sqrt{\Delta}}{a}$

$$\phi_y\left\{\psi_{xx} + \frac{2b}{a}\psi_{xy} + \frac{c}{a}\psi_{yy}\right\} = \psi_y\left\{\phi_{xx} + \frac{2b}{a}\phi_{xy} + \frac{c}{a}\phi_{yy}\right\}$$

Infine si moltiplica per *a*

$$\phi_y\{a\psi_{xx} + 2b\psi_{xy} + c\psi_{yy}\} = \psi_y\{a\phi_{xx} + 2b\phi_{xy} + c\phi_{yy}\}$$

e si conclude che

$$[\mathcal{L}^-, \mathcal{L}^+] = 0 \Leftrightarrow \mathcal{L}^+[\Lambda^-] = \mathcal{L}^-[\Lambda^+] \Leftrightarrow \phi_y\mathcal{L}[\psi] = \psi_y\mathcal{L}[\phi]$$

dove $\mathcal{L}$ è l'operatore della parte principale della PDE ed $\mathcal{L}^-, \mathcal{L}^+$ sono la sua fattorizzazione; in effetti dividendo per $a\phi_y\psi_y$, in virtù delle uguaglianze $\mathcal{L}[\phi] = a\mathcal{L}^-[\Lambda^+]\phi_y$, $\mathcal{L}[\psi] = a\mathcal{L}^+[\Lambda^-]\psi_y$ si trova proprio $\mathcal{L}^+[\Lambda^-] = \mathcal{L}^-[\Lambda^+]$. Nello specifico si ha che

$$\mathcal{L}^-[\Lambda^+] = 0 \Leftrightarrow \psi_y(\phi_{xx}\phi_y - \phi_{xy}\phi_x) + \psi_x(\phi_{yy}\phi_x - \phi_{xy}\phi_y) = 0$$

$$\mathcal{L}^+[\Lambda^-] = 0 \Leftrightarrow \phi_y(\psi_{xx}\psi_y - \psi_{xy}\psi_x) + \phi_x(\psi_{yy}\psi_x - \psi_{xy}\psi_y) = 0$$

Si osserva che le condizioni sono verificate sicuramente quando i coefficienti di (*) sono costanti, essendo nulle le derivate seconde degli invarianti. Dividendo entrambe le espressioni per $\phi_y\psi_y$ e procedendo come prima

$$\mathcal{L}^-[\Lambda^+] = 0 \Leftrightarrow \mathcal{L}[\phi] = 0$$

$$\mathcal{L}^+[\Lambda^-] = 0 \Leftrightarrow \mathcal{L}[\psi] = 0$$

Qualora valessero contemporaneamente si avrebbe correttamente che $\phi_y\mathcal{L}[\psi] = \psi_y\mathcal{L}[\phi]$ e

$$\mathcal{L} = a\mathcal{L}^-\mathcal{L}^+ = a\mathcal{L}^+\mathcal{L}^- \Leftrightarrow \mathcal{L}^+[\Lambda^-] = \mathcal{L}^-[\Lambda^+] = 0 \Leftrightarrow \mathcal{L}[\phi] = \mathcal{L}[\psi] = 0$$

*Es.4*

*Nell'esempio 3 si avevano gli invarianti*

$$\phi(x,y) = \ln(e^x - 1) - 2e^{-\frac{y}{2}} \qquad \psi(x,y) = \ln(e^x - 1) - 2e^{-\frac{y}{2}} - x$$

*La condizione di commutatività dei rispettivi operatori*

$$\phi_y \{ a\psi_{xx} + 2b\psi_{xy} + c\psi_{yy}\} = \psi_y\{ a\phi_{xx} + 2b\phi_{xy} + c\phi_{yy}\}$$

*è facilmente verificata, essendo*

$$\phi_x = \frac{e^x}{e^x - 1} \qquad \psi_x = \frac{1}{e^x - 1} \qquad \phi_y = \psi_y = e^{-y/2}$$

$$\phi_{xx} = \psi_{xx} = -\frac{e^x}{(e^x - 1)^2} \qquad \phi_{xy} = \psi_{xy} = 0 \qquad \phi_{yy} = \psi_{yy} = -\frac{1}{2} e^{-\frac{y}{2}}$$

*Tuttavia* $\mathcal{L}[\phi] = \mathcal{L}[\psi] = -\left(2 + e^{\frac{y}{2}}\right) \neq 0$ *e infatti si aveva*

$$\mathcal{L}^+[\Lambda^-] = \mathcal{L}^-[\Lambda^+] = -\frac{1}{8} e^{\frac{y}{2}} csch^2\left(\frac{x}{2}\right)\left(2 + e^{\frac{y}{2}}\right) \neq 0$$

*Es.5*

$$u_{tt} + 4t u_{tx} + 3t^2 u_{xx} = 0$$

*Si tenga a mente che a = 1. Da* $(\Lambda^\pm)^2 + 4t\Lambda^\pm + 3t^2 = 0$ *si ha che* $\Lambda^+ = -t$, $\Lambda^- = -3t$, *per cui*

$$\mathcal{L}^+ = \partial_t + t\partial_x \qquad \mathcal{L}^- = \partial_t + 3t\partial_x$$

*e da* $\frac{dy}{dx} = -\Lambda^\pm$ *si trovano i due invarianti e quindi la trasformazione:*

$$\phi(t,x) = x - \frac{t^2}{2} = \xi \qquad\qquad \phi^{-1}(\xi,\eta) = (\xi - \eta)^{\frac{1}{2}} = t$$

$$\psi(t,x) = x - \frac{3t^2}{2} = \eta \qquad\qquad \psi^{-1}(\xi,\eta) = \frac{3\xi - \eta}{2} = x$$

*Da* $\mathcal{L}^-[\Lambda^+] = -1 \neq \mathcal{L}^+[\Lambda^-] = -3$ *segue che* $[\mathcal{L}^-, \mathcal{L}^+] \neq 0$. *Si calcolano ancora:* $\phi_x = \psi_x = 1$ *e*

$\mathcal{L}^-[\phi] = 2t \quad \mathcal{L}^+[\psi] = -2t \quad$ *da cui* $\quad \mathcal{L}^-[\phi]\mathcal{L}^+[\psi] = -4t^2 = -4(\xi - \eta)$. *Non resta che applicare*

$$\mathcal{L} = \Phi^{-1} \circ \{ a [ (\mathcal{L}^-[\phi]\mathcal{L}^+[\psi])\partial_{\xi\eta} + (\mathcal{L}^-[\Lambda^+]\phi_x)\partial_\xi + (\mathcal{L}^+[\Lambda^-]\psi_x)\partial_\eta ]\}$$

$$\mathcal{L} = -4(\xi - \eta)\partial_{\xi\eta} - \partial_\xi - 3\partial_\eta$$

*Pertanto la forma canonica dell'equazione sarà*

$$4(\xi - \eta)U_{\xi\eta} + U_\xi + 3U_\eta = 0$$

*Benché ridurre in forma canonica la PDE con questo metodo sia molto sbrigativo, non lo è altrettanto trovare un integrale generale, in quanto l'equazione, a prima vista molto semplice, non si riduce ad una ODE. Ci si può accontentare di una soluzione particolare non banale che si ottiene nel seguente modo. Si scrive la funzione incognita come:*

$$U(\xi, \eta) = F(\xi, \eta) + G(\xi) + H(\eta)$$

*in modo da avere*

$$U_\xi = F_\xi + G'(\xi) \qquad U_\eta = F_\eta + H'(\eta) \qquad U_{\xi\eta} = F_{\xi\eta}$$

*Andando a sostituire nell'equazione*

$$4(\xi - \eta)F_{\xi\eta} + F_\xi + G'(\xi) + 3F_\eta + 3H'(\eta) = 0$$

$$\left(4\xi F_{\xi\eta} + 3F_\eta\right) - \left(4\eta F_{\xi\eta} - F_\xi\right) = -G'(\xi) - 3H'(\eta)$$

*A questo punto si impone*

$$4\eta (F_\xi)_\eta - F_\xi = 3H'(\eta)$$

$$4\xi (F_\eta)_\xi + 3F_\eta = -G'(\xi)$$

*Le due equazioni differenziali ordinarie, lineari e del prim'ordine hanno soluzione:*

$$F_\xi = \eta^{\frac{1}{4}} c_1(\xi) + \frac{3}{4} \eta^{\frac{1}{4}} \int \eta^{-\frac{5}{4}} H'(\eta)\, d\eta$$

$$F_\eta = \xi^{-\frac{3}{4}} c_2(\eta) - \frac{1}{4} \xi^{-\frac{3}{4}} \int \xi^{-\frac{1}{4}} G'(\xi)\, d\xi$$

*Integrando ancora una volta si ottengono due espressioni per $F(\xi, \eta)$ che dovranno risultare equivalenti*

$$F(\xi, \eta) = \eta^{\frac{1}{4}} C_1(\xi) + \frac{3}{4} \xi \eta^{\frac{1}{4}} \int \eta^{-\frac{5}{4}} H'(\eta)\, d\eta$$

$$F(\xi, \eta) = \xi^{-\frac{3}{4}} C_2(\eta) - \frac{1}{4} \xi^{-\frac{3}{4}} \eta \int \xi^{-\frac{1}{4}} G'(\xi)\, d\xi$$

*Uguagliando i primi termini*

$$\eta^{\frac{1}{4}} C_1(\xi) = \xi^{-\frac{3}{4}} C_2(\eta)$$

*si ha che* $C_1(\xi) = c_0 \xi^{-\frac{3}{4}}$, $C_2(\eta) = c_0 \eta^{\frac{1}{4}}$ *con,* $c_0 = cost.$, *ovvero*

$$F(\xi, \eta) = c_0 \left(\frac{\eta}{\xi^3}\right)^{\frac{1}{4}} + \cdots$$

*Dovrà ancora essere*

$$\frac{3}{4}\xi \left(\eta^{\frac{1}{4}} \int \eta^{-\frac{5}{4}} H'(\eta)\, d\eta\right) = -\frac{1}{4}\eta \left(\xi^{-\frac{3}{4}} \int \xi^{-\frac{1}{4}} G'(\xi)\, d\xi\right)$$

$$\xi^{-\frac{3}{4}} \int \xi^{-\frac{1}{4}} G'(\xi)\, d\xi = k_1 \xi \quad \Rightarrow \quad \int \xi^{-\frac{1}{4}} G'(\xi)\, d\xi = k_1 \xi^{\frac{7}{4}}$$

$$\eta^{\frac{1}{4}} \int \eta^{-\frac{5}{4}} H'(\eta)\, d\eta = k_2 \eta \quad \Rightarrow \quad \int \eta^{-\frac{5}{4}} H'(\eta)\, d\eta = k_2\, \eta^{\frac{3}{4}}$$

*da cui*

$$\frac{3}{4}\xi\, k_2 \eta = -\frac{1}{4} k_1 \xi\, \eta \quad \Rightarrow \quad k_1 = -3k_2$$

$$\int \xi^{-\frac{1}{4}} G'(\xi)\, d\xi = -3k_2 \xi^{\frac{7}{4}}$$

*Derivando ambo i membri e poi integrando si ricava*

$$G'(\xi) = -\frac{21}{4}k_2\xi \implies G(\xi) = -\frac{21}{8}k_2\xi^2$$

*Si ripete il procedimento per il secondo integrale*

$$H'(\eta) = \frac{3}{4}k_2\eta \implies H(\eta) = \frac{3}{8}k_2\eta^2$$

*Si ponga per comodità $k_2 = \frac{8}{3}c_1$*

$$F(\xi,\eta) = c_0\left(\frac{\eta}{\xi^3}\right)^{\frac{1}{4}} + 2c_1\xi\eta \qquad G(\xi) = -7c_1\xi^2 \qquad H(\eta) = c_1\eta^2$$

*Sommando i termini si ottiene infine*

$$U(\xi,\eta) = F(\xi,\eta) + G(\xi) + H(\eta) = c_0\left(\frac{\eta}{\xi^3}\right)^{\frac{1}{4}} + c_1(\eta^2 + 2\xi\eta - 7\xi^2)$$

*che può anche scriversi*

$$U(\xi,\eta) = c_0\left(\frac{\eta}{\xi^3}\right)^{\frac{1}{4}} + c_1\left[(\xi-\eta)^2 - 2\left(\frac{3\xi-\eta}{2}\right)^2\right]$$

*Inoltre, poiché $u_1 = k_1 t$ e $u_2 = k_2 x$ sono soluzioni particolari e ricordando le trasformazioni inverse per tornare alle variabili iniziali, si arriva ad una soluzione particolare del tipo*

$$u(t,x) = a + bt + cx + d(t^4 - 2x^2) + e\left[\frac{3t^2 - 2x}{(t^2 - 2x)^3}\right]^{\frac{1}{4}}$$

*dove a,b,c,d,e sono costanti arbitrarie.*

Nel **problema parabolico**

$$\Delta = b^2 - ac = 0 \implies \Lambda = -\frac{b}{a}$$

la parte principale viene fattorizzata nel prodotto di due operatori identici

$$\mathcal{L}^2[u] = a\mathcal{L}\mathcal{L}[u] = au_{xx} + 2bu_{xy} + cu_{yy}$$

Di conseguenza si ha una sola famiglia di caratteristiche, soluzione di:

$$\mathcal{L}[\phi] = 0$$

Si sceglie ancora $\psi$ in modo da completare la mappa di transizione e tale per cui:

$$\mathcal{L}[\psi] \neq 0 \quad \text{e} \quad \mathcal{L}^2[\psi] = 0$$

ovvero

$$a\psi_{xx} + 2b\psi_{xy} + c\psi_{yy} = 0$$

Nel caso parabolico non si pone il problema della commutatività degli operatori in quanto si ha solo un operatore che commuta con se stesso; comunemente al caso iperbolico resta il problema se il prodotto dell'operatore per se stesso sia effettivamente uguale alla forma quadratica nelle derivate parziali, che altro non è che la parte principale di (*), senza termini aggiuntivi di grado inferiore. Avendo definito

$$\mathcal{L} = \left(\frac{\partial}{\partial x} - \Lambda \frac{\partial}{\partial y}\right)$$

applicando l'operatore a se stesso:

$$\mathcal{L}\mathcal{L} = a^{-1}\mathcal{L}^2 - \mathcal{L}[\Lambda]\partial_y$$

da cui:

$$\mathcal{L}^2[u] = a\mathcal{L}\mathcal{L}[u] \Leftrightarrow \mathcal{L}[\Lambda] = 0$$

L'osservazione più ovvia che si può fare è che se $\mathcal{L}[\Lambda] = 0$, $\Lambda(x,y)$ è invariante rispetto a $\mathcal{L}$ e come tale si può esprimere come funzione di $\phi(x,y)$ e viceversa purché $\Lambda \neq cost.$:

$$\mathcal{L}[\phi] = \mathcal{L}[f(\Lambda)] = f'(\Lambda)\mathcal{L}[\Lambda] = 0$$

$$\mathcal{L}[\phi] = (\partial_x - \Lambda\,\partial_y)\phi = \phi_x - f^{-1}(\phi)\,\phi_y = 0 \Rightarrow f^{-1}(\phi) = \frac{\phi_x}{\phi_y} = \Lambda \Rightarrow \phi = f(\Lambda)$$

Inoltre $\mathcal{L}[\Lambda]$ si esprime rispetto agli invarianti ad $\mathcal{L}$ ed $\mathcal{L}^2$ come

$$\mathcal{L}[\Lambda] = \frac{1}{\phi_y^3}\{\phi_x^2\phi_{yy} - 2\phi_x\phi_y\phi_{xy} + \phi_y^2\phi_{xx}\}$$

quindi

$$\mathcal{L}[\Lambda] = 0 \Leftrightarrow \phi_x^2\phi_{yy} - 2\phi_x\phi_y\phi_{xy} + \phi_y^2\phi_{xx} = 0$$

condizione sicuramente verificata quando le derivate seconde di $\phi$ sono nulle ($\Lambda = cost.$). Dividendo per $\phi_y^2$, sapendo che $\frac{\phi_x}{\phi_y} = \Lambda$,

$$\Lambda^2\phi_{yy} - 2\Lambda\phi_{xy} + \phi_{xx} = 0$$

ma $\Lambda = -\frac{b}{a}$ e $c = \frac{b^2}{a}$

$$\phi_{xx} + \frac{2b}{a}\phi_{xy} + \frac{c}{a}\phi_{yy} = 0$$

Si conclude moltiplicando per $a$

$$a\phi_{xx} + 2b\phi_{xy} + c\phi_{yy} = 0$$

ovvero

$$\mathcal{L}^2 = a\mathcal{L}\mathcal{L} \Leftrightarrow \mathcal{L}[\Lambda] = 0 \Leftrightarrow \mathcal{L}^2[\phi] = 0$$

E' abbastanza immediato: essendo $\mathcal{L}[\Lambda] = 0$, $\mathcal{L}^2 = a\mathcal{L}\mathcal{L}$, per cui se $\phi$ è invariante rispetto a $\mathcal{L}$ lo sarà anche rispetto a $\mathcal{L}^2 = a\mathcal{L}\mathcal{L}$.

*Es.1*

$$\mathcal{L}^2 = x^2 \partial_{xx} + 2xy\, \partial_{xy} + y^2 \partial_{yy} = x^2 \mathcal{L}\mathcal{L} \quad con \quad \mathcal{L} = \partial_x + \frac{y}{x}\partial_y$$

*Siccome $\mathcal{L}^2 = a\mathcal{L}\mathcal{L}$, $\mathcal{L}[\Lambda] = 0$ per cui si può scegliere come invariante una funzione di $\Lambda = -\frac{y}{x}$, ad esempio*

$$\phi(x,y) = \frac{y}{x}$$

*Volendo verificare che*

$$\mathcal{L}^2[\phi] = x^2 \phi_{xx} + 2xy\phi_{xy} + y^2 \phi_{yy} = 0$$

*si calcolano*

$$\phi_{xx} = \frac{2y}{x^3} \qquad \phi_{xy} = -\frac{1}{x^2} \qquad \phi_{yy} = 0$$

*da cui sostituendo*

$$\mathcal{L}^2[\phi] = x^2 \frac{2y}{x^3} + 2xy \left(-\frac{1}{x^2}\right) + y^2 \cdot 0 = \frac{2y}{x} - \frac{2y}{x} = 0$$

*E' chiaro che qualsiasi funzione di $\phi(x,y) = \frac{y}{x}$ è candidata a soluzione particolare della PDE.*

Considerando:

$$\mathcal{L}(\mathcal{L}[\psi]) = a^{-1}\mathcal{L}^2[\psi] - \mathcal{L}[\Lambda]\psi_y$$

$$\mathcal{L}(\mathcal{L}[\psi]) = \frac{1}{\phi_y^3}\left\{\phi_x^2(\psi_{yy}\phi_y - \phi_{yy}\psi_y) - 2\phi_x\phi_y(\psi_{xy}\phi_y - \phi_{xy}\psi_y) + \phi_y^2(\psi_{xx}\phi_y - \phi_{xx}\psi_y)\right\}$$

si ha che la condizione $\mathcal{L}^2[\psi] = 0$ si traduce in

$$\mathcal{L}(\mathcal{L}[\psi]) + \mathcal{L}[\Lambda]\psi_y = 0$$

Qualora $\psi$ fosse funzione della sola $x$ si avrebbe che $\mathcal{L}(\mathcal{L}[\psi]) = \psi_{xx} = 0$, ovvero la scelta deve ricadere su una funzione lineare in $x$ per essere compatibile con i requisiti del problema. Lo stesso vale nel caso in cui $\psi = \psi(y)$, infatti

$$\mathcal{L}(\mathcal{L}[\psi]) = \frac{\phi_x^2}{\phi_y^2}\psi_{yy} - 2\frac{\phi_x}{\phi_y}\psi_{xy} + \psi_{xx} - \frac{1}{\phi_y^3}\left\{\phi_x^2\phi_{yy} - 2\phi_x\phi_y\phi_{xy} + \phi_y^2\phi_{xx}\right\}\psi_y$$
$$= a^{-1}\mathcal{L}^2[\psi] - \mathcal{L}[\Lambda]\psi_y = -\mathcal{L}[\Lambda]\psi_y$$

purché $\psi_{yy} = 0$, come si sarebbe dedotto facilmente anche da $a\psi_{xx} + 2b\psi_{xy} + c\psi_{yy} = 0$. E' chiaro che se $\mathcal{L}^2 = a\mathcal{L}\mathcal{L}$, cioè $\mathcal{L}[\Lambda] = 0$, allora $\mathcal{L}(\mathcal{L}[\psi]) = 0$. Riassumendo, è possibile avere $\mathcal{L}^2[\psi] = a\mathcal{L}(\mathcal{L}[\psi]) = 0$ se $\psi = \psi(x)$ pur essendo in generale $\mathcal{L}^2 \neq a\mathcal{L}\mathcal{L}$, oppure per $\mathcal{L}[\Lambda] = 0$.

Dopo queste osservazioni si può procedere al calcolo dell'operatore parabolico nelle nuove coordinate. Si applichi dunque l'operatore differenziale su $U(\xi, \eta) = U(\phi(x,y), \psi(x,y))$:

$$\mathcal{L}[U] = \frac{\partial U}{\partial \xi} \mathcal{L}[\phi] + \frac{\partial U}{\partial \eta} \mathcal{L}[\psi] = \frac{\partial U}{\partial \eta} \mathcal{L}[\psi]$$

essendo $\phi$ costante rispetto a $\mathcal{L}$.

Segue che

$$\frac{\partial}{\partial \eta} = \frac{\mathcal{L}}{\mathcal{L}[\psi]}$$

$$\frac{\partial^2}{\partial \eta^2} = \frac{\mathcal{L}}{\mathcal{L}[\psi]}\left[\frac{\mathcal{L}}{\mathcal{L}[\psi]}\right] = \frac{\mathcal{L}\mathcal{L}}{(\mathcal{L}[\psi])^2} - \frac{\mathcal{L}(\mathcal{L}[\psi])}{(\mathcal{L}[\psi])^2}\frac{\mathcal{L}}{\mathcal{L}[\psi]} = \frac{\mathcal{L}\mathcal{L}}{(\mathcal{L}[\psi])^2} - \frac{\mathcal{L}(\mathcal{L}[\psi])}{(\mathcal{L}[\psi])^2}\frac{\partial}{\partial \eta}$$

Supponendo per il momento di essere nella condizione in cui $\mathcal{L}(\mathcal{L}[\psi]) = 0$ si arriva a

$$\mathcal{L}\mathcal{L} = (\mathcal{L}[\psi])^2 \frac{\partial^2}{\partial \eta^2} \qquad (^{\circ\circ})$$

Tenendo presente che $\Lambda = -\frac{b}{a}$, $\mathcal{L} = \left(\frac{\partial}{\partial x} + \frac{b}{a}\frac{\partial}{\partial y}\right) = \frac{1}{a}\left(a\frac{\partial}{\partial x} + b\frac{\partial}{\partial y}\right)$, $\mathcal{L}[\psi] = \frac{1}{a}\left(a\psi_x + b\psi_y\right)$, $b^2 = ac$, ma anche che $\Lambda = \frac{\phi_x}{\phi_y}$ e dunque $\mathcal{L}[\psi] = \frac{1}{\phi_y}\left(\phi_y\psi_x - \phi_x\psi_y\right) = -\frac{1}{\phi_y}\left|\frac{\partial(\phi,\psi)}{\partial(x,y)}\right|$, si può scrivere:

$$(\mathcal{L}[\psi])^2 = a^{-1}\left(a\psi_x^2 + 2b\psi_x\psi_y + c\psi_y^2\right) = (\phi_y)^{-2}\left|\frac{\partial(\phi,\psi)}{\partial(x,y)}\right|^2$$

Se l'operatore della parte principale è proprio uguale (a meno del coefficiente *a*) al quadrato dell'operatore $\mathcal{L}$ ($\mathcal{L}[\Lambda] = 0$, condizione certamente valida per (*) a coefficienti costanti)

$$\mathcal{L}^2[u] = a\mathcal{L}\mathcal{L}[u] = \left(a\psi_x^2 + 2b\psi_x\psi_y + c\psi_y^2\right)\frac{\partial^2 U}{\partial \eta^2}$$

tenendo conto dei termini di grado inferiore:

$$\left(a\psi_x^2 + 2b\psi_x\psi_y + c\psi_y^2\right)U_{\eta\eta} + \mathcal{F}(\xi, \eta, U_\xi, U_\eta, U) = 0$$

da cui la forma canonica:

$$U_{\eta\eta} = F(\xi, \eta, U_\xi, U_\eta, U)$$

Questa volta si ha la corrispondenza:

$$\mathcal{LL}[u(x,y)] \xrightarrow{(\phi,\psi)} \partial_{\eta\eta} U(\xi,\eta)$$

*Es. Caso Omogeneo*

$$U_{\eta\eta} = 0 \Rightarrow U_\eta = F(\xi) \Rightarrow U(\xi,\eta) = \eta F(\xi) + G(\xi)$$

*Rientra in questa casistica l'esempio 1, con $\phi(x,y) = \frac{y}{x} = \xi$, si completa la mappa con $\psi(x,y) = x = \eta$, $\mathcal{L}[\psi] = 1$, dalla (°°) $\mathcal{LL} = \partial_{\eta\eta}$, $\mathcal{L}^2 = x^2 \mathcal{LL} = \eta^2 \partial_{\eta\eta}$. L'integrale generale nelle coordinate $(x,y)$ sarà*

$$u(x,y) = x\, F\left(\frac{y}{x}\right) + G\left(\frac{y}{x}\right)$$

*Es.2*

$$\mathcal{L}^2 = y^2 \partial_{xx} - 2y\, \partial_{xy} + \partial_{yy} = \mathcal{LL} + \partial_x \;,\; \mathcal{L} = \partial_y - y\partial_x$$

*Data ora la seguente equazione differenziale*

$$\mathcal{L}^2[u] = y^2 u_{xx} - 2y\, u_{xy} + u_{yy} = u_x + 6y$$

*Si può sfruttare la precedente relazione:*

$$\tilde{\mathcal{L}}^2[u] = \mathcal{LL}[u] = y^2 u_{xx} - 2y\, u_{xy} + u_{yy} - u_x = 6y$$

*Ora si può applicare la (°°):*

$$\tilde{\mathcal{L}}^2[u] = U_{\eta\eta}\, (\mathcal{L}[\psi])^2 = 6y$$

*Risolvendo $\mathcal{L}[\phi] = \phi_y - y\phi_x = 0$ si trova $\phi(x,y) = x + \frac{y^2}{2} = \xi$ e completiamo la mappa con $x = \xi - \frac{y^2}{2}$, $y = \eta$, da cui $\mathcal{L}[\psi] = \partial_y y - y\partial_x y = 1$ ; pertanto la forma canonica sarà:*

$$U_{\eta\eta} = 6\eta$$

*L'integrale generale è semplice da calcolare:*

$$U(\xi,\eta) = \eta^3 + \eta\, F(\xi) + G(\xi)$$

*Tornando alle coordinate iniziali:*

$$u(x,y) = y^3 + y\, F\left(x + \frac{y^2}{2}\right) + G\left(x + \frac{y^2}{2}\right)$$

Riprendendo i conti e tendendo a mente le uguaglianze viste in precedenza si ha

$$\frac{\partial^2}{\partial \eta^2} = \frac{\mathcal{LL}}{(\mathcal{L}[\psi])^2} - \frac{\mathcal{L}(\mathcal{L}[\psi])}{(\mathcal{L}[\psi])^2} \frac{\partial}{\partial \eta} = \frac{\mathcal{LL}}{(\mathcal{L}[\psi])^2} + \frac{\mathcal{L}[\Lambda]\psi_y}{(\mathcal{L}[\psi])^2} \frac{\partial}{\partial \eta}$$

ovvero

$$\mathcal{LL} = a^{-1} \mathcal{L}^2 - \mathcal{L}[\Lambda]\partial_y = (\mathcal{L}[\psi])^2\, \partial_{\eta\eta} - \mathcal{L}[\Lambda]\psi_y\, \partial_\eta$$

$$a^{-1}\mathcal{L}^2 = (\mathcal{L}[\psi])^2 \partial_{\eta\eta} - \mathcal{L}[\Lambda]\psi_y \partial_\eta + \mathcal{L}[\Lambda]\partial_y$$

$$a^{-1}\mathcal{L}^2 = (\mathcal{L}[\psi])^2 \partial_{\eta\eta} - \mathcal{L}[\Lambda]\psi_y \partial_\eta + \mathcal{L}[\Lambda]\phi_y \partial_\xi + \mathcal{L}[\Lambda]\psi_y \partial_\eta$$

da cui segue l'espressione per *l'operatore parabolico della parte principale nelle variabili* $(\xi, \eta)$

$$\mathcal{L}^2 = \Phi^{-1} \circ \left\{ a \left[ (\mathcal{L}[\psi])^2 \frac{\partial^2}{\partial \eta^2} + (\mathcal{L}[\Lambda]\phi_y) \frac{\partial}{\partial \xi} \right] \right\}$$

La bontà dell'espressione è confermata dal fatto che la presenza delle derivate al suo interno, ossia la sua forma, non dipende dalla scelta dell'arbitraria funzione $\psi$ (si ricordi che $\mathcal{L}[\psi] \neq 0$). In aggiunta, poiché

$$\mathcal{L}(\mathcal{L}[\phi]) = a^{-1}\mathcal{L}^2[\phi] - \mathcal{L}[\Lambda]\phi_y = 0 \quad \Longrightarrow \quad \mathcal{L}^2[\phi] = a\mathcal{L}[\Lambda]\phi_y$$

si può anche scrivere

$$\mathcal{L}^2 = \Phi^{-1} \circ \left[ a(\mathcal{L}[\psi])^2 \frac{\partial^2}{\partial \eta^2} + \mathcal{L}^2[\phi] \frac{\partial}{\partial \xi} \right]$$

*Es. Equazione del calore (o di diffusione)*

$$u_t = k\, u_{xx} \qquad k > 0 \ \ cost.$$

*E' evidente che la PDE di tipo parabolico si trovi già nella sua forma canonica.*

*Es 3.*

$$xy^3 u_{xx} - 2x^2 y^2\, u_{xy} + x^3 y\, u_{yy} = y^3 u_x + x^3 u_y$$

*Innanzitutto si evidenzia la parte principale*

$$y^2 u_{xx} - 2xy\, u_{xy} + x^2 u_{yy} = \frac{1}{xy}\left( y^3 u_x + x^3 u_y \right)$$

*Il problema è parabolico poiché*

$$y^2 \Lambda^2 - 2xy\, \Lambda + x^2 = (y\Lambda - x)^2 = 0 \quad \Longrightarrow \quad \Lambda = \frac{x}{y}$$

*Si trova la famiglia di caratteristiche risolvendo $\frac{dy}{dx} = -\Lambda = -\frac{x}{y}$ e si sceglie ad esempio $\psi = x$ in modo da completare la mappa*

$$\phi(x,y) = x^2 + y^2 \quad = \xi \qquad\qquad \phi^{-1}(\xi,\eta) = \eta \qquad\qquad = x$$
$$\psi(x,y) = x \qquad\quad = \eta \qquad\qquad \psi^{-1}(\xi,\eta) = (\xi - \eta^2)^{\frac{1}{2}} = y$$

*Si calcoli ancora:* $\mathcal{L}[\psi] = 1, \phi_x = 2x, \phi_y = 2y, \psi_x = 1, , \psi_y = 0,$ *e*

$$\mathcal{L}[\Lambda] = \left(\partial_x - \frac{x}{y}\partial_y\right)\frac{x}{y} = \frac{x^2 + y^2}{y^3} \neq 0$$

*Dall'espressione dell'operatore parabolico segue che*

$$\mathcal{L}^2 = a(\mathcal{L}[\psi])^2 \partial_{\eta\eta} + a\mathcal{L}[\Lambda]\phi_y \partial_\xi = y^2\partial_{\eta\eta} + y^2\frac{x^2+y^2}{y^3}2y\,\partial_\xi = (\xi - \eta^2)\partial_{\eta\eta} + 2\xi\partial_\xi$$

*oppure, sapendo che* $\phi_{xx} = \phi_{yy} = 2$, $\phi_{xy} = 0$ *e quindi* $\mathcal{L}^2[\phi] = 2(x^2 + y^2)$, *si può usare la seconda formulazione*

$$\mathcal{L}^2 = a(\mathcal{L}[\psi])^2 \partial_{\eta\eta} + \mathcal{L}^2[\phi]\partial_\xi = y^2\partial_{\eta\eta} + 2(x^2 + y^2) = (\xi - \eta^2)\partial_{\eta\eta} + 2\xi\partial_\xi$$

*Il secondo membro dell'equazione,* $\frac{1}{xy}\left(y^3\partial_x + x^3\partial_y\right)u(x,y)$, *si trasforma mediante*

$$\partial_x = \phi_x \partial_\xi + \psi_x \partial_\eta = 2x\partial_\xi + \partial_\eta$$

$$\partial_y = \phi_y \partial_\xi + \psi_y \partial_\xi = 2y\partial_\xi$$

$$\frac{1}{xy}\left(y^3\partial_x + x^3\partial_y\right) = \frac{1}{xy}\left(2xy^3\partial_\xi + y^3\partial_\eta + 2yx^3\partial_\xi\right) = 2(x^2 + y^2)\partial_\xi + \frac{y^2}{x}\partial_\eta = 2\xi\partial_\xi + \frac{\xi - \eta^2}{\eta}\partial_\eta$$

*Applicando gli operatori alla funzione incognita nelle nuove variabili si ha*

$$(\xi - \eta^2)U_{\eta\eta} + 2\xi U_\xi = 2\xi U_\xi + \frac{\xi - \eta^2}{\eta}U_\eta$$

*da cui la forma canonica*

$$U_{\eta\eta} = \frac{U_\eta}{\eta}$$

*E' facile integrare la ODE che ha come soluzione*

$$U(\xi,\eta) = \eta^2 F(\xi) + G(\xi)$$

*Si conclude tornando alle variabili* $(x,y)$

$$u(x,y) = x^2 F(x^2 + y^2) + G(x^2 + y^2)$$

## *Alcune osservazioni sulla mappa inversa* $\Phi^{-1}$

Si consideri la trasformazione di coordinate $\Phi(\phi(x,y), \psi(x,y))$

$$\begin{cases} \xi = \phi(x,y) \\ \eta = \psi(x,y) \end{cases}$$

Differenziando si ha

$$d\xi = d\phi(x,y) = \phi_x\, dx + \phi_y\, dy$$
$$d\eta = d\psi(x,y) = \psi_x\, dx + \psi_y\, dy$$

Poiché $d\xi$ e $d\eta$ sono versori dello *spazio cotangente* (*duale* dello spazio tangente) è possibile moltiplicarli secondo le regole del *prodotto wedge* (prodotto esterno)

$$d\xi \wedge d\xi = d\eta \wedge d\eta = 0 \qquad d\xi \wedge d\eta = -\, d\eta \wedge d\xi$$

Allora

$$d\xi \wedge d\eta = (\phi_x\, dx + \phi_y\, dy) \wedge (\psi_x\, dx + \psi_y\, dy) =$$

$$= \phi_x \psi_x\, dx \wedge dx + \phi_x \psi_y\, dx \wedge dy + \phi_y \psi_x\, dy \wedge dx + \phi_y \psi_y\, dy \wedge dy =$$

$$= (\phi_x \psi_y - \phi_y \psi_x) dx \wedge dy = \begin{vmatrix} \phi_x & \phi_y \\ \psi_x & \psi_y \end{vmatrix} dx \wedge dy = |J_\Phi(x,y)|\, dx \wedge dy$$

dove $|J_\Phi(x,y)| = \left|\frac{\partial(\phi,\psi)}{\partial(x,y)}\right|$ è il determinante jacobiano della trasformazione $\Phi$.

Si consideri ora la trasformazione inversa $\Phi^{-1} = (\phi^{-1}(\xi,\eta), \psi^{-1}(\xi,\eta))$

$$\begin{cases} x = \phi^{-1}(\xi,\eta) \\ y = \psi^{-1}(\xi,\eta) \end{cases}$$

Differenziando si ha

$$dx = d\phi^{-1}(\xi,\eta) = \phi^{-1}_\xi\, d\xi + \phi^{-1}_\eta\, d\eta$$
$$dy = d\psi^{-1}(\xi,\eta) = \psi^{-1}_\xi\, d\xi + \psi^{-1}_\eta\, d\eta$$

Questa volta si avrà che

$$dx \wedge dy = \begin{vmatrix} \phi^{-1}_\xi & \phi^{-1}_\eta \\ \psi^{-1}_\xi & \psi^{-1}_\eta \end{vmatrix} d\xi \wedge d\eta = |J_{\Phi^{-1}}(\xi,\eta)|\, d\xi \wedge d\eta$$

dove $|J_{\Phi^{-1}}(\xi,\eta)| = \left|\frac{\partial(\phi^{-1},\psi^{-1})}{\partial(\xi,\eta)}\right|$ è il determinante jacobiano della trasformazione inversa $\Phi^{-1}$.

Invertendo la prima relazione

$$dx \wedge dy = \frac{1}{|J_\Phi(x,y)|} d\xi \wedge d\eta$$

e confrontandola con la seconda

$$dx \wedge dy = |J_{\Phi^{-1}}(\xi,\eta)| d\xi \wedge d\eta$$

si vede che

$$|J_{\Phi^{-1}}(\xi,\eta)| = \frac{1}{|J_\Phi(x,y)|}$$

Siccome per qualsiasi matrice invertibile $\boldsymbol{A}$ vale la regola

$$|\boldsymbol{A^{-1}}| = \frac{1}{|\boldsymbol{A}|}$$

si conclude che

$$\boldsymbol{J}_{\Phi^{-1}}(\xi,\eta) = [\boldsymbol{J}_\Phi(x,y)]^{-1}$$

ovvero: la matrice jacobiana della trasformazione inversa $\Phi^{-1}$ è l'inversa della matrice jacobiana della trasformazione diretta $\Phi$ ; vale ovviamente il viceversa

$$\boldsymbol{J}_\Phi(x,y) = [\boldsymbol{J}_{\Phi^{-1}}(\xi,\eta)]^{-1}$$

Si stabiliscono quindi delle relazioni che saranno utili nel seguito:

$$\phi_x = \frac{\psi_\eta^{-1}}{|J_{\Phi^{-1}}(\xi,\eta)|} \qquad \phi_y = -\frac{\phi_\eta^{-1}}{|J_{\Phi^{-1}}(\xi,\eta)|}$$

$$\psi_x = -\frac{\psi_\xi^{-1}}{|J_{\Phi^{-1}}(\xi,\eta)|} \qquad \psi_y = \frac{\phi_\xi^{-1}}{|J_{\Phi^{-1}}(\xi,\eta)|}$$

$$\phi_\xi^{-1} = \frac{\psi_y}{|J_\Phi(x,y)|} \qquad \phi_\eta^{-1} = -\frac{\phi_y}{|J_\Phi(x,y)|}$$

$$\psi_\xi^{-1} = -\frac{\psi_x}{|J_\Phi(x,y)|} \qquad \psi_\eta^{-1} = \frac{\phi_x}{|J_\Phi(x,y)|}$$

*Problema Iperbolico*

Si ricordi che

$$\mathcal{L}^- = \mathcal{L}^-[\phi]\, \partial_\xi \qquad \mathcal{L}^+ = \mathcal{L}^+[\psi]\, \partial_\eta$$

Mediante un calcolo esplicito e utilizzando le relazioni precedenti

$$\mathcal{L}^-[\phi] = \phi_x - \Lambda^- \phi_y = \phi_x - \frac{\psi_x}{\psi_y}\phi_y = \frac{1}{\psi_y}(\phi_x \psi_y - \psi_x \phi_y) = \frac{|J_\Phi(x,y)|}{\psi_y} = \frac{1}{\phi_\xi^{-1}}$$

$$\mathcal{L}^+[\psi] = \psi_x - \Lambda^+ \psi_y = \psi_x - \frac{\phi_x}{\phi_y}\psi_y = \frac{1}{\phi_y}(\psi_x \phi_y - \phi_x \psi_y) = -\frac{|J_\Phi(x,y)|}{\phi_y} = \frac{1}{\phi_\eta^{-1}}$$

$$\mathcal{L}^- = \frac{1}{\phi_\xi^{-1}}\,\partial_\xi \qquad \mathcal{L}^+ = \frac{1}{\phi_\eta^{-1}}\partial_\eta$$

Il prodotto degli operatori è

$$\mathcal{L}^-\mathcal{L}^+ = \frac{1}{\phi_\xi^{-1}}\,\partial_\xi\left(\frac{1}{\phi_\eta^{-1}}\,\partial_\eta\right) = \frac{1}{\phi_\xi^{-1}}\left(-\frac{\phi_{\xi\eta}^{-1}}{(\phi_\eta^{-1})^2}\,\partial_\eta + \frac{1}{\phi_\eta^{-1}}\,\partial_{\xi\eta}\right) = \frac{1}{\phi_\xi^{-1}\phi_\eta^{-1}}\left(\partial_{\xi\eta} - \phi_{\xi\eta}^{-1}\,\mathcal{L}^+\right)$$

$$\mathcal{L}^+\mathcal{L}^- = \frac{1}{\phi_\eta^{-1}}\,\partial_\eta\left(\frac{1}{\phi_\xi^{-1}}\,\partial_\xi\right) = \frac{1}{\phi_\eta^{-1}}\left(-\frac{\phi_{\xi\eta}^{-1}}{(\phi_\xi^{-1})^2}\,\partial_\xi + \frac{1}{\phi_\xi^{-1}}\,\partial_{\xi\eta}\right) = \frac{1}{\phi_\xi^{-1}\phi_\eta^{-1}}\left(\partial_{\xi\eta} - \phi_{\xi\eta}^{-1}\,\mathcal{L}^-\right)$$

Si può calcolare il commutatore

$$[\mathcal{L}^-,\mathcal{L}^+] = \mathcal{L}^-\mathcal{L}^+ - \mathcal{L}^+\mathcal{L}^- = \frac{\phi_{\xi\eta}^{-1}}{\phi_\xi^{-1}\phi_\eta^{-1}}\,(\mathcal{L}^- - \mathcal{L}^+)$$

Affinché si annulli, esclusa la soluzione $\mathcal{L}^- = \mathcal{L}^+$ che non è accettabile, dovrà essere $\phi_{\xi\eta}^{-1} = 0$. Si possono calcolare per verifica $\mathcal{L}^-[\Lambda^+]$ e $\mathcal{L}^+[\Lambda^-]$, sapendo che

$$\Lambda^+ = \frac{\phi_x}{\phi_y} = -\frac{\psi_\eta^{-1}}{\phi_\eta^{-1}} \qquad\qquad \Lambda^- = \frac{\psi_x}{\psi_y} = -\frac{\psi_\xi^{-1}}{\phi_\xi^{-1}}$$

$$\mathcal{L}^-[\Lambda^+] = \frac{1}{\phi_\xi^{-1}}\,\partial_\xi\left(-\frac{\psi_\eta^{-1}}{\phi_\eta^{-1}}\right) = \frac{\phi_{\xi\eta}^{-1}\psi_\eta^{-1} - \psi_{\xi\eta}^{-1}\phi_\eta^{-1}}{\phi_\xi^{-1}(\phi_\eta^{-1})^2} = \frac{\phi_\xi^{-1}\psi_\eta^{-1}\phi_{\xi\eta}^{-1} - \phi_\xi^{-1}\phi_\eta^{-1}\psi_{\xi\eta}^{-1}}{(\phi_\xi^{-1}\phi_\eta^{-1})^2}$$

$$\mathcal{L}^+[\Lambda^-] = \frac{1}{\phi_\eta^{-1}}\,\partial_\eta\left(-\frac{\psi_\xi^{-1}}{\phi_\xi^{-1}}\right) = \frac{\phi_{\xi\eta}^{-1}\psi_\xi^{-1} - \psi_{\xi\eta}^{-1}\phi_\xi^{-1}}{(\phi_\xi^{-1})^2\phi_\eta^{-1}} = \frac{\phi_\eta^{-1}\psi_\xi^{-1}\phi_{\xi\eta}^{-1} - \phi_\xi^{-1}\phi_\eta^{-1}\psi_{\xi\eta}^{-1}}{(\phi_\xi^{-1}\phi_\eta^{-1})^2}$$

Dato che

$$[\mathcal{L}^-,\mathcal{L}^+] = 0 \iff \mathcal{L}^-[\Lambda^+] = \mathcal{L}^+[\Lambda^-]$$

$$\phi_\xi^{-1}\psi_\eta^{-1}\phi_{\xi\eta}^{-1} - \phi_\xi^{-1}\phi_\eta^{-1}\psi_{\xi\eta}^{-1} = \phi_\eta^{-1}\psi_\xi^{-1}\phi_{\xi\eta}^{-1} - \phi_\xi^{-1}\phi_\eta^{-1}\psi_{\xi\eta}^{-1}$$

$$\phi_\xi^{-1}\psi_\eta^{-1}\phi_{\xi\eta}^{-1} = \phi_\eta^{-1}\psi_\xi^{-1}\phi_{\xi\eta}^{-1} \implies -\frac{\psi_\eta^{-1}}{\phi_\eta^{-1}}\,\phi_{\xi\eta}^{-1} = -\frac{\psi_\xi^{-1}}{\phi_\xi^{-1}}\,\phi_{\xi\eta}^{-1}$$

$$\Lambda^+\phi_{\xi\eta}^{-1} = \Lambda^-\phi_{\xi\eta}^{-1}$$

Poiché $\Lambda^+ \neq \Lambda^-$, si conferma che

$$[\mathcal{L}^-,\mathcal{L}^+] = 0 \iff \mathcal{L}^-[\Lambda^+] = \mathcal{L}^+[\Lambda^-] \iff \phi_{\xi\eta}^{-1} = 0$$

$\phi^{-1}_{\xi\eta} = 0$ significa che in generale $\phi^{-1}(\xi,\eta) = f(\xi) + g(\eta)$, ovvero $x = f(\phi(x,y)) + g(\psi(x,y))$ da cui

$$\psi(x,y) = g^{-1}\left(x - f(\phi(x,y))\right)$$

Si scelga un operatore $\mathcal{L}^+$ il quale ammetterà come invariante $\phi(x,y)$, oppure lo si costruisca a partire da quest'ultimo ($\mathcal{L}^+ = \partial_x + \Lambda^+ \partial_y$, $\Lambda^+ = \frac{\phi_x}{\phi_y}$); qualsiasi operatore $\mathcal{L}^-$ che abbia come invariante $\psi(x,y) = f\big(g(\phi(x,y)) - x\big)$, dove $f, g$ sono funzioni arbitrariamente scelte, commuta con $\mathcal{L}^+$. Quindi $\mathcal{L}^-$ avrà la forma $\mathcal{L}^- = \partial_x + \Lambda^- \partial_y$, con

$$\Lambda^- = \frac{\psi_x}{\psi_y} = \frac{f'(g(\phi) - x)\,(g'(\phi)\phi_x - 1)}{f'(g(\phi) - x)g'(\phi)\phi_y} = \frac{g'(\phi)\phi_x - 1}{g'(\phi)\phi_y}$$

$\Lambda^-$ non dipende da $f$: d'altronde è sufficiente che sia invariante il suo argomento, per cui per costruire $\mathcal{L}^-$ basta prendere $\psi = g(\phi) - x$. Manipolando l'espressione si ottiene la relazione

$$\Lambda^- = \left(1 - \frac{1}{g'(\phi)\phi_x}\right)\Lambda^+$$

Tornando alle espressioni per $\mathcal{L}^-[\Lambda^+]$ e $\mathcal{L}^+[\Lambda^-]$, posto $\phi^{-1}_{\xi\eta} \neq 0$, affinché sia $\mathcal{L}^-[\Lambda^+] = 0$ ($\mathcal{L}^+[\Lambda^-] \neq 0$), $\phi^{-1}_{\xi\eta}\psi^{-1}_\eta = \psi^{-1}_{\xi\eta}\phi^{-1}_\eta$ e dunque

$$\frac{\psi^{-1}_{\xi\eta}}{\psi^{-1}_\eta} = \frac{\phi^{-1}_{\xi\eta}}{\phi^{-1}_\eta} \quad \Longrightarrow \quad \frac{\partial}{\partial \xi}\ln\left(\frac{\psi^{-1}_\eta}{\phi^{-1}_\eta}\right) = 0 \quad \Longrightarrow \quad \frac{\psi^{-1}_{\xi\eta}}{\phi^{-1}_{\xi\eta}} = \frac{\psi^{-1}_\eta}{\phi^{-1}_\eta} = -\Lambda^+ = f(\eta)$$

ovvero $\Lambda^+ = f(\psi(x,y))$, per cui $\mathcal{L}^-[\Lambda^+] = \mathcal{L}^-[f(\psi)] = f'(\psi)\mathcal{L}^-[\psi] = 0$.

Se $\mathcal{L}^+[\Lambda^-] = 0$ ($\mathcal{L}^-[\Lambda^+] \neq 0$) invece, $\phi^{-1}_{\xi\eta}\psi^{-1}_\xi = \psi^{-1}_{\xi\eta}\phi^{-1}_\xi$, quindi

$$\frac{\psi^{-1}_{\xi\eta}}{\psi^{-1}_\xi} = \frac{\phi^{-1}_{\xi\eta}}{\phi^{-1}_\xi} \quad \Longrightarrow \quad \frac{\partial}{\partial \eta}\ln\left(\frac{\psi^{-1}_\xi}{\phi^{-1}_\xi}\right) = 0 \quad \Longrightarrow \quad \frac{\psi^{-1}_{\xi\eta}}{\phi^{-1}_{\xi\eta}} = \frac{\psi^{-1}_\xi}{\phi^{-1}_\xi} = -\Lambda^- = g(\xi)$$

ossia $\Lambda^- = g(\phi(x,y))$, infatti $\mathcal{L}^+[\Lambda^-] = \mathcal{L}^+[g(\phi)] = g'(\phi)\mathcal{L}^+[\phi] = 0$.

E' chiaro che date le ipotesi, le quali ci consentono di dividere per una quantità non nulla, non è possibile che sia $\mathcal{L}^-[\Lambda^+] = \mathcal{L}^+[\Lambda^-] = 0$, per cui si avrebbe $\phi^{-1}_{\xi\eta} = 0$ e $\frac{\psi^{-1}_{\xi\eta}}{\phi^{-1}_{\xi\eta}} = -\Lambda^+ = -\Lambda^-$.

Si supponga ora che $\phi^{-1}_{\xi\eta} = 0$; affinché $\mathcal{L}^-[\Lambda^+]$ e $\mathcal{L}^+[\Lambda^-]$ siano entrambi nulli devono valere contemporaneamente le condizioni

$$\psi^{-1}_{\xi\eta}\phi^{-1}_\eta = 0 \qquad \psi^{-1}_{\xi\eta}\phi^{-1}_\xi = 0$$

$\phi^{-1}_\xi = \phi^{-1}_\eta = 0$ non ha senso e si conclude che

$$\mathcal{L} = a\,\mathcal{L}^-\mathcal{L}^+ = a\,\mathcal{L}^+\mathcal{L}^- \iff \mathcal{L}^-[\Lambda^+] = \mathcal{L}^+[\Lambda^-] = 0 \iff \phi^{-1}_{\xi\eta} = \psi^{-1}_{\xi\eta} = 0$$

*Es. 1*

*Nell'esempio 3 del caso iperbolico si aveva*

$$\phi(x,y) = \ln(e^x - 1) - 2e^{-\frac{y}{2}} = \xi \qquad \phi^{-1}(\xi,\eta) = \xi - \eta = x$$

$$\psi(x,y) = \ln(e^x - 1) - 2e^{-\frac{y}{2}} - x = \eta \qquad \psi^{-1}(\xi,\eta) = -2\ln\left\{\tfrac{1}{2}\left[\ln(e^{\xi-\eta}-1) - \xi\right]\right\} = y$$

*Si osserva che $\psi(x,y) = \phi(x,y) - x$ per cui $g(\phi) = \phi$, $g'(\phi) = 1$. Inoltre*

$$\Lambda^- = \left(1 - \frac{1}{g'(\phi)\phi_x}\right)\Lambda^+ = \left(1 - \frac{e^x-1}{e^x}\right)\frac{e^{x+\frac{y}{2}}}{e^x-1} = e^{-x}\frac{e^{x+\frac{y}{2}}}{e^x-1} = \frac{e^{\frac{y}{2}}}{e^x-1}$$

*A partire dalla trasformazione di coordinate inversa si ha che $\phi^{-1}_{\xi\eta} = 0$, condizione che garantisce la commutatività degli operatori $\mathcal{L}^-, \mathcal{L}^+$, mentre, anche senza un calcolo esplicito, $\psi^{-1}_{\xi\eta} \neq 0$ dal momento che non è nella forma $\psi^{-1}(\xi,\eta) = f(\xi) + g(\eta)$; difatti ci si trova nel caso in cui $\mathcal{L}^-[\Lambda^+] = \mathcal{L}^+[\Lambda^-] \neq 0$.*

*Es. 2*

*Si ha l'equazione, sulla falsariga dell'esempio 5 del problema iperbolico*

$$t^2 u_{tt} + 4tx u_{tx} + 3x^2 u_{xx} = 0$$

*Risolvendo $t^2(\Lambda^\pm)^2 + 4tx\Lambda^\pm + 3x^2 = 0$ si trovano $\Lambda^+ = -\frac{x}{t}$, $\Lambda^- = -\frac{3x}{t}$ e quindi da $\frac{dx}{dt} = -\Lambda^\pm$ si hanno gli invarianti*

$$\phi(x,y) = \frac{x}{t} = \xi \qquad \phi^{-1}(\xi,\eta) = \left(\frac{\xi}{\eta}\right)^{\frac{1}{2}} = t$$

$$\psi(x,y) = \frac{x}{t^3} = \eta \qquad \psi^{-1}(\xi,\eta) = \left(\frac{\xi^3}{\eta}\right)^{\frac{1}{2}} = x$$

*L'intenzione è, a parte risolvere l'equazione, quella di verificare le uguaglianze viste in precedenza. Quindi, a partire dalla mappa inversa $\Phi^{-1}$, si calcolano*

$$\phi^{-1}_\xi = \frac{1}{2}\xi^{-\frac{1}{2}}\eta^{-\frac{1}{2}} \qquad \phi^{-1}_\eta = -\frac{1}{2}\xi^{\frac{1}{2}}\eta^{-\frac{3}{2}} \qquad \phi^{-1}_{\xi\eta} = -\frac{1}{4}\xi^{-\frac{1}{2}}\eta^{-\frac{3}{2}}$$

$$\psi^{-1}_\xi = \frac{3}{2}\xi^{\frac{1}{2}}\eta^{-\frac{1}{2}} \qquad \psi^{-1}_\eta = -\frac{1}{2}\xi^{\frac{3}{2}}\eta^{-\frac{3}{2}} \qquad \psi^{-1}_{\xi\eta} = -\frac{3}{4}\xi^{\frac{1}{2}}\eta^{-\frac{3}{2}}$$

*Si ha effettivamente che*

$$\Lambda^+ = -\frac{\psi^{-1}_\eta}{\phi^{-1}_\eta} = -\xi = -\frac{x}{t} \qquad \Lambda^- = -\frac{\psi^{-1}_\xi}{\phi^{-1}_\xi} = -3\xi = -\frac{3x}{t}$$

*Si controlli ora il rapporto*

$$\frac{\psi^{-1}_{\xi\eta}}{\phi^{-1}_{\xi\eta}} = -\frac{3}{4}\xi^{\frac{1}{2}}\eta^{-\frac{3}{2}}\left(-4\xi^{\frac{1}{2}}\eta^{\frac{3}{2}}\right) = 3\xi = \frac{\psi^{-1}_\xi}{\phi^{-1}_\xi} = -\Lambda^-$$

*Ciò significa che* $\mathcal{L} = a\,\mathcal{L}^+\mathcal{L}^- \neq a\,\mathcal{L}^-\mathcal{L}^+ \Leftrightarrow \mathcal{L}^-[\Lambda^+] \neq \mathcal{L}^+[\Lambda^-] = 0$ *e difatti* $\mathcal{L}^-[\Lambda^+] = -\frac{2x}{t^2} \neq 0$. *Si può procedere con*

$$\mathcal{L}^+\mathcal{L}^- = \frac{1}{\phi_\xi^{-1}\phi_\eta^{-1}}\left(\partial_{\xi\eta} - \phi_{\xi\eta}^{-1}\,\mathcal{L}^-\right)$$

*per cui servono*

$$\mathcal{L}^- = \frac{1}{\phi_\xi^{-1}}\partial_\xi = 2\xi^{\frac{1}{2}}\eta^{\frac{1}{2}}\partial_\xi \qquad \phi_{\xi\eta}^{-1}\,\mathcal{L}^- = -\frac{1}{4}\xi^{-\frac{1}{2}}\eta^{-\frac{3}{2}}\left(2\xi^{\frac{1}{2}}\eta^{\frac{1}{2}}\right)\partial_\xi = -\frac{1}{2\eta}\partial_\xi$$

$$\frac{1}{\phi_\xi^{-1}\phi_\eta^{-1}} = 2\xi^{\frac{1}{2}}\eta^{\frac{1}{2}}\left(-2\xi^{-\frac{1}{2}}\eta^{\frac{3}{2}}\right) = -4\eta^2$$

*Sostituendo*

$$\mathcal{L}^+\mathcal{L}^- = -4\eta^2\left(\partial_{\xi\eta} + \frac{1}{2\eta}\partial_\xi\right) = -2\eta\left(2\eta\,\partial_{\xi\eta} + \partial_\xi\right)$$

*Poiché* $a = t^2 = \frac{\xi}{\eta}$ *si conclude che*

$$\mathcal{L} = -2\xi\left(2\eta\,\partial_{\xi\eta} + \partial_\xi\right)$$

*Applicando l'operatore alla funzione* $U(\xi,\eta)$ *si ha la forma canonica*

$$2\eta U_{\xi\eta} + U_\xi = 0$$

*Se i conti sono giusti, mediante l'espressione*

$$\mathcal{L} = \Phi^{-1}\circ\left\{a\left[\left(\mathcal{L}^-[\phi]\,\mathcal{L}^+[\psi]\right)\partial_{\xi\eta} + \left(\mathcal{L}^-[\Lambda^+]\,\phi_x\right)\partial_\xi + \left(\mathcal{L}^+[\Lambda^-]\,\psi_x\right)\partial_\eta\right]\right\}$$

*si arriva allo stesso risultato. Servono:* $\mathcal{L}^-[\phi] = \frac{2x}{t^2}$, $\mathcal{L}^+[\psi] = -\frac{2x}{t^4}$, $\phi_x = \frac{1}{t}$, $\psi_x = \frac{1}{t^3}$, *per cui*

$$\mathcal{L} = \Phi^{-1}\circ\left\{t^2\left[-\frac{4x^2}{t^6}\partial_{\xi\eta} - \frac{2x}{t^3}\partial_\xi\right]\right\} = \Phi^{-1}\circ\left\{\left[-\frac{4x}{t}\cdot\frac{x}{t^3}\partial_{\xi\eta} - \frac{2x}{t}\partial_\xi\right]\right\} = -4\xi\eta\partial_{\xi\eta} - 2\xi\partial_\xi$$

*Per integrare l'equazione si scrive*

$$2\eta U_{\xi\eta} + U_\xi = \partial_\xi\left(2\eta U_\eta + U\right) = 0 \quad\Rightarrow\quad U_\eta + \frac{1}{2\eta}U = g(\eta)$$

*L'ODE lineare del prim'ordine si risolve facilmente con il metodo di Lagrange*

$$U(\xi,\eta) = \eta^{-\frac{1}{2}}\left[\int \eta^{\frac{1}{2}}\,g(\eta)d\eta + f(\xi)\right]$$

$$U(\xi,\eta) = \eta^{-\frac{1}{2}}F(\xi) + G(\eta)$$

*da cui applicando* $\Phi$

$$u(t,x) = \left(\frac{t^3}{x}\right)^{\frac{1}{2}}F\left(\frac{x}{t}\right) + G\left(\frac{x}{t^3}\right)$$

*Si noti come* $u(t,x) = \frac{x}{t^3} = \psi(x,y)$ *sia una soluzione particolare, in accordo con la regola*

$$\mathcal{L}^+[\Lambda^-] = 0 \iff \mathcal{L}[\psi] = 0$$

*Un altro modo di procedere consiste nel fattorizzare*

$$2\eta U_{\xi\eta} + U_\xi = [(2\eta\partial_\eta + 1)\partial_\xi] U(\xi,\eta) = (\mathcal{L}_1 \mathcal{L}_2) U(\xi,\eta)$$

*Poiché gli operatori commutano si cerca una soluzione nella forma*

$$U = U_1 + U_2 \quad con \quad \mathcal{L}_1[U_1] = 0 \ e \ \mathcal{L}_2[U_2] = 0$$

*Si ha infatti*

$$(\mathcal{L}_1 \mathcal{L}_2) U = (\mathcal{L}_1 \mathcal{L}_2)(U_1 + U_2) = (\mathcal{L}_1 \mathcal{L}_2)U_1 + (\mathcal{L}_1 \mathcal{L}_2)U_2 = \mathcal{L}_2(\mathcal{L}_1 [U_1]) + \mathcal{L}_1(\mathcal{L}_2 [U_2]) = 0$$

*Da* $\mathcal{L}_2[U_2] = \partial_\xi U_2 = 0$ *si ha subito* $U_2 = g(\eta)$, *mentre*

$$\mathcal{L}_1[U_1] = (2\eta\partial_\eta + 1)U_1 = 2\eta\,(U_1)_\eta + U_1 = 0$$

*Si tratta di risolvere un'ODE a variabili separabili*

$$\int \frac{dU_1}{U_1} = -\frac{1}{2}\int \frac{d\eta}{\eta} \Rightarrow \ln U_1 = -\frac{1}{2}\ln\eta + f(\xi) \Rightarrow \ln U_1 \eta^{\frac{1}{2}} = f(\xi) \Rightarrow U_1 = \eta^{-\frac{1}{2}} F(\xi)$$

*Sommando si trova proprio*

$$U(\xi,\eta) = U_1 + U_2 = \eta^{-\frac{1}{2}} F(\xi) + G(\eta)$$

<u>*Problema Parabolico*</u>

Si comincia calcolando

$$\mathcal{L}[\psi] = \psi_x - \Lambda\psi_y = \psi_x - \frac{\phi_x}{\phi_y}\psi_y = \frac{1}{\phi_y}(\psi_x \phi_y - \phi_x \psi_y) = -\frac{|J_\Phi(x,y)|}{\phi_y} = \frac{1}{\phi_\eta^{-1}}$$

per cui

$$\mathcal{L} = \mathcal{L}[\psi]\partial_\eta = \frac{1}{\phi_\eta^{-1}}\,\partial_\eta$$

$$\mathcal{L}\mathcal{L} = \frac{1}{\phi_\eta^{-1}}\,\partial_\eta\left(\frac{1}{\phi_\eta^{-1}}\,\partial_\eta\right) = \frac{1}{\phi_\eta^{-1}}\left(-\frac{\phi_{\eta\eta}^{-1}}{\phi_\eta^{-1}}\partial_\eta + \frac{1}{\phi_\eta^{-1}}\partial_{\eta\eta}\right) = \frac{1}{(\phi_\eta^{-1})^2}(\partial_{\eta\eta} - \phi_{\eta\eta}^{-1}\mathcal{L})$$

Pertanto

$$\mathcal{L}(\mathcal{L}[\psi]) = \frac{1}{(\phi_\eta^{-1})^2}(\partial_{\eta\eta} - \phi_{\eta\eta}^{-1}\mathcal{L})\,\eta = -\frac{\phi_{\eta\eta}^{-1}}{(\phi_\eta^{-1})^3}$$

La condizione che portava alla (°°) sarà

$$\mathcal{L}(\mathcal{L}[\psi]) = 0 \iff \phi_{\eta\eta}^{-1} = 0$$

Siccome

$$\Lambda = -\frac{\psi_\eta^{-1}}{\phi_\eta^{-1}}$$

allora

$$\mathcal{L}[\Lambda] = \frac{1}{\phi_\eta^{-1}}\, \partial_\eta \left(-\frac{\psi_\eta^{-1}}{\phi_\eta^{-1}}\right) = \frac{\phi_{\eta\eta}^{-1}\,\psi_\eta^{-1} - \psi_{\eta\eta}^{-1}\,\phi_\eta^{-1}}{\left(\phi_\eta^{-1}\right)^3}$$

Ricordando che $\mathcal{L}^2[\psi] = 0$ implicava $\mathcal{L}(\mathcal{L}[\psi]) + \mathcal{L}[\Lambda]\psi_y = 0$

$$\frac{\phi_{\eta\eta}^{-1}}{\left(\phi_\eta^{-1}\right)^3} = \frac{\phi_{\eta\eta}^{-1}\,\psi_\eta^{-1} - \psi_{\eta\eta}^{-1}\,\phi_\eta^{-1}}{\left(\phi_\eta^{-1}\right)^3} \cdot \frac{\phi_\xi^{-1}}{|J_{\Phi^{-1}}(\xi,\eta)|}$$

$$\phi_\xi^{-1}\psi_\eta^{-1}\phi_{\eta\eta}^{-1} - \phi_\eta^{-1}\psi_\xi^{-1}\phi_{\eta\eta}^{-1} = \phi_\xi^{-1}\psi_\eta^{-1}\phi_{\eta\eta}^{-1} - \phi_\xi^{-1}\phi_\eta^{-1}\psi_{\eta\eta}^{-1}$$

da cui la proprietà

$$\mathcal{L}^2[\psi] = 0 \iff \psi_\xi^{-1}\phi_{\eta\eta}^{-1} = \phi_\xi^{-1}\psi_{\eta\eta}^{-1}$$

Qualora $\phi_{\eta\eta}^{-1} = 0$ ($\mathcal{L}(\mathcal{L}[\psi]) = 0$), $\phi_\xi^{-1}\psi_{\eta\eta}^{-1} = 0$, per cui si potrebbe avere $\psi_{\eta\eta}^{-1} = 0$ (che data l'uguaglianza $\mathcal{L}(\mathcal{L}[\psi]) + \mathcal{L}[\Lambda]\psi_y = 0$ corrisponderà alla condizione $\mathcal{L}[\Lambda] = 0$), oppure se $\psi_\xi^{-1} \neq 0$ ($\psi_x \neq 0$, non può essere $\psi_x = \psi_y = 0$), $\phi_\xi^{-1} = 0$ o ancora, $\phi_\xi^{-1} = \psi_{\eta\eta}^{-1} = 0$. Se $\psi_\xi^{-1} \neq 0$ ($\psi_x \neq 0$) e $\phi_{\eta\eta}^{-1} = 0$ ($\mathcal{L}(\mathcal{L}[\psi]) = 0$), allora $\phi_\xi^{-1}\psi_{\eta\eta}^{-1} = 0$. Supponendo $\psi_{\eta\eta}^{-1} \neq 0$ ($\mathcal{L}[\Lambda] \neq 0$), non resta che porre $\phi_\xi^{-1} = 0$ ($\psi_y = 0$) ovvero $\psi = \psi(x)$. Da $\phi_{\eta\eta}^{-1} = 0$ segue che $\phi^{-1}(\xi,\eta) = \eta f(\xi) + g(\xi)$, ma $\phi_\xi^{-1} = \eta f'(\xi) + g'(\xi) = 0$ e dunque $f'(\xi) = g'(\xi) = 0$, cioè $\phi^{-1}(\xi,\eta) = c_1'\eta + c_2' = x$, con $c_1', c_2' = cost$. Invertendo si trova $\eta = c_1 x + c_2 = \psi(x)$: la mappa $\Phi$ è stata completata con una funzione lineare della sola variabile $x$, o meglio, $\mathcal{L}(\mathcal{L}[\psi]) = 0$ implica, nel caso in cui $\mathcal{L}[\Lambda] \neq 0$, che si è scelta una funzione del tipo (si osservi che $\psi_y = 0$, in accordo con $\mathcal{L}(\mathcal{L}[\psi]) + \mathcal{L}[\Lambda]\psi_y = 0$). Similmente, se $\psi_{\eta\eta}^{-1} = 0$, $\psi_\xi^{-1}\phi_{\eta\eta}^{-1} = 0$; segue che se $\phi_{\eta\eta}^{-1} = 0$, $\mathcal{L}[\Lambda] = 0$ e di conseguenza $\mathcal{L}(\mathcal{L}[\psi]) = a^{-1}\mathcal{L}^2[\psi] = 0$. Se $\phi_{\eta\eta}^{-1} \neq 0$, $\mathcal{L}[\Lambda] \neq 0$, dovrà essere $\psi_\xi^{-1} = 0$ ($\psi_x = 0$) implicando $\phi_\xi^{-1} \neq 0$ ($\psi_y \neq 0$); si potrà infine avere, posto ancora $\phi_\xi^{-1} \neq 0$, $\psi_\xi^{-1} = \phi_{\eta\eta}^{-1} = 0$. Supponendo ad esempio di trovarsi in questo caso ($\mathcal{L}[\Lambda] = 0$), con $\phi_\xi^{-1} \neq 0$ e $\psi_\xi^{-1} = \phi_{\eta\eta}^{-1} = \psi_{\eta\eta}^{-1} = 0$, si avrà che $\phi_{\eta\eta}^{-1} = 0 \Rightarrow \phi^{-1}(\xi,\eta) = \eta f(\xi) + g(\xi)$, ma con $\phi_\xi^{-1} = \eta f'(\xi) + g'(\xi) \neq 0$, per cui almeno uno tra $f(\xi), g(\xi)$ dovrà non essere un numero. Inoltre $\psi_\xi^{-1} = 0 \Rightarrow \psi^{-1}(\xi,\eta) = h(\eta)$ da cui derivando $\psi_{\eta\eta}^{-1} = h''(\eta) = 0$ e poi integrando $h'(\eta) = c_1' \Rightarrow h(\eta) = c_1'\eta + c_2' = \psi^{-1}(\eta) = y$, si conclude che $\eta = c_1 y + c_2 = \psi(y)$; qualora $\psi = \psi(y)$ e $\mathcal{L}[\Lambda] = 0$, $\psi(y)$ deve essere una funzione lineare nella sua variabile affinché si abbia $\mathcal{L}(\mathcal{L}[\psi]) = a^{-1}\mathcal{L}^2[\psi] = 0$. E' chiaro che si sarebbe giunti alla stessa conclusione se $\phi_{\eta\eta}^{-1}$ fosse stato non nullo ($\mathcal{L}(\mathcal{L}[\psi]) \neq 0$ e quindi $\mathcal{L}[\Lambda] \neq 0$ essendo $\psi_y \neq 0$), in quanto deve comunque valere la condizione $\mathcal{L}^2[\psi] = 0 \Rightarrow \psi''(y) = 0$. I risultati dedotti partendo dalla mappa inversa coincidono con la precedente trattazione sul caso parabolico.

Ponendo uguale a zero l'espressione per $\mathcal{L}[\Lambda]$ si ottiene l'uguaglianza $\phi^{-1}_{\eta\eta}\psi^{-1}_\eta = \psi^{-1}_{\eta\eta}\phi^{-1}_\eta$, da cui

$$\frac{\psi^{-1}_{\eta\eta}}{\psi^{-1}_\eta} = \frac{\phi^{-1}_{\eta\eta}}{\phi^{-1}_\eta} \quad\Rightarrow\quad \frac{\partial}{\partial\eta}\ln\left(\frac{\psi^{-1}_\eta}{\phi^{-1}_\eta}\right) = 0 \quad\Rightarrow\quad \frac{\psi^{-1}_\eta}{\phi^{-1}_\eta} = -\Lambda = f(\xi)$$

ovvero $\Lambda = f(\phi(x,y))$ affinché sia $\mathcal{L}[\Lambda] = \mathcal{L}[f(\phi)] = f'(\phi)\mathcal{L}[\phi] = 0$, come si è già visto. Esplicitando nella precedente

$$\phi^{-1}_{\eta\eta} = \frac{\phi^{-1}_\eta}{\psi^{-1}_\eta}\,\psi^{-1}_{\eta\eta}$$

ed inserendo in $\phi^{-1}_\xi \psi^{-1}_{\eta\eta} = \psi^{-1}_\xi \phi^{-1}_{\eta\eta}$, che deve necessariamente essere soddisfatta, si ha

$$\phi^{-1}_\xi \psi^{-1}_{\eta\eta} = \psi^{-1}_\xi \frac{\phi^{-1}_\eta}{\psi^{-1}_\eta}\,\psi^{-1}_{\eta\eta}$$

Svolgendo i calcoli

$$\phi^{-1}_\xi \psi^{-1}_\eta \psi^{-1}_{\eta\eta} = \psi^{-1}_\xi \phi^{-1}_\eta \psi^{-1}_{\eta\eta} \quad\Rightarrow\quad (\phi^{-1}_\xi \psi^{-1}_\eta - \psi^{-1}_\xi \phi^{-1}_\eta)\psi^{-1}_{\eta\eta} = |J_{\Phi^{-1}}(\xi,\eta)|\,\psi^{-1}_{\eta\eta} = 0$$

Siccome la trasformazione è invertibile, $|J_{\Phi^{-1}}(\xi,\eta)| \neq 0$, da cui la soluzione $\psi^{-1}_{\eta\eta} = 0$; ma se $\psi^{-1}_{\eta\eta}$ è nullo lo sarà anche $\phi^{-1}_{\eta\eta} = \frac{\phi^{-1}_\eta}{\psi^{-1}_\eta}\psi^{-1}_{\eta\eta}$ e si completa così la dimostrazione di quanto dedotto precedentemente per altre vie:

$$\mathcal{L}^2 = a\mathcal{L}\mathcal{L} \;\Leftrightarrow\; \mathcal{L}[\Lambda] = 0 \;\Leftrightarrow\; \phi^{-1}_{\eta\eta} = \psi^{-1}_{\eta\eta} = 0$$

*Es.1*

*Nell'esempio 3 del caso parabolico si aveva*

$$\begin{aligned}\phi(x,y) &= x^2 + y^2 &&= \xi & \phi^{-1}(\xi,\eta) &= \eta &&= x\\ \psi(x,y) &= x &&= \eta & \psi^{-1}(\xi,\eta) &= (\xi - \eta^2)^{\frac{1}{2}} &&= y\end{aligned}$$

*Si calcolino*

$$\phi^{-1}_\xi = 0 \qquad \phi^{-1}_\eta = 1 \qquad \phi^{-1}_{\eta\eta} = 0$$

$$\psi^{-1}_\xi = \frac{1}{2}(\xi - \eta^2)^{-\frac{1}{2}} \qquad \psi^{-1}_\eta = -\eta(\xi - \eta^2)^{-\frac{1}{2}} \qquad \psi^{-1}_{\eta\eta} = -\xi(\xi - \eta^2)^{-\frac{3}{2}}$$

$$\Lambda = -\frac{\psi^{-1}_\eta}{\phi^{-1}_\eta} = \frac{\eta}{(\xi - \eta^2)^{\frac{1}{2}}} = \frac{x}{y}$$

*Da $\phi^{-1}_{\eta\eta} = 0$ segue che $\mathcal{L}(\mathcal{L}[\psi]) = 0$, difatti $\mathcal{L}[\psi] = \frac{1}{\phi^{-1}_\eta} = 1$, $\mathcal{L}[1] = 0$; tuttavia $\psi^{-1}_{\eta\eta} \neq 0$, per cui $\mathcal{L}[\Lambda] \neq 0$, e dunque deve essere $\phi^{-1}_\xi = 0$ ($\phi_y = 0$). Deve comunque valere l'uguaglianza ($\mathcal{L}^2[\psi] = 0$)*

$$\psi_\xi^{-1} \phi_{\eta\eta}^{-1} = \phi_\xi^{-1} \psi_{\eta\eta}^{-1} \quad \Rightarrow \quad \frac{1}{2}(\xi - \eta^2)^{-\frac{1}{2}} \cdot 0 = 0 \cdot \left[-\xi(\xi - \eta^2)^{-\frac{3}{2}}\right]$$

*Poiché $\phi_{\eta\eta}^{-1} = 0$ l'espressione per $\mathcal{L}[\Lambda]$ si riduce a*

$$\mathcal{L}[\Lambda] = -\frac{\psi_{\eta\eta}^{-1}}{\left(\phi_\eta^{-1}\right)^2} = \frac{\xi}{(\xi - \eta^2)^{\frac{3}{2}}} = \frac{x^2 + y^2}{y^3}$$

*Sapendo che $a = y^2 = \xi - \eta^2$ e*

$$\phi_y = -\frac{\phi_\eta^{-1}}{|J_{\Phi^{-1}}(\xi, \eta)|} = -\frac{1}{-\frac{1}{2}(\xi - \eta^2)^{-\frac{1}{2}}} = 2(\xi - \eta^2)^{\frac{1}{2}} = 2y$$

*si può applicare*

$$\mathcal{L}^2 = \Phi^{-1} \circ \left\{ a \left[ (\mathcal{L}[\psi])^2 \partial_{\eta\eta} + (\mathcal{L}[\Lambda] \phi_y) \partial_\xi \right] \right\}$$

$$\mathcal{L}^2 = (\xi - \eta^2) \left[ \partial_{\eta\eta} + \left( \frac{\xi}{(\xi - \eta^2)^{\frac{3}{2}}} \cdot 2(\xi - \eta^2)^{\frac{1}{2}} \right) \partial_\xi \right] = (\xi - \eta^2) \partial_{\eta\eta} + 2\xi \partial_\xi$$

*Procedendo dalla trasformazione inversa $\Phi^{-1}$, oltre ad aver verificato il funzionamento dei teoremi visti, si è giunti ad un'espressione dell'operatore $\mathcal{L}^2$ nelle coordinate $(\xi, \eta)$ che coincide proprio con quella trovata nell'esempio 3. Il problema poteva risolversi anche scegliendo*

$$\begin{aligned} \phi(x,y) &= x^2 + y^2 = \xi & \phi^{-1}(\xi, \eta) &= (\xi - \eta^2)^{\frac{1}{2}} = x \\ \psi(x,y) &= y = \eta & \psi^{-1}(\xi, \eta) &= \eta = y \end{aligned}$$

*Essendo $\psi_y = 1$ ed avendo già calcolato $\mathcal{L}[\Lambda]$ è immediato conoscere*

$$\mathcal{L}(\mathcal{L}[\psi]) = -\mathcal{L}[\Lambda]\psi_y = -\frac{x^2 + y^2}{y^3} = -\frac{\xi}{\eta^3} \neq 0$$

*difatti stavolta $\phi_{\eta\eta}^{-1} \neq 0$; a partire da $\phi^{-1}(\xi, \eta)$ si ha proprio che*

$$\mathcal{L}(\mathcal{L}[\psi]) = -\frac{\phi_{\eta\eta}^{-1}}{\left(\phi_\eta^{-1}\right)^3} = -\left[ -\frac{\xi}{(\xi - \eta^2)^{\frac{3}{2}}} \left( -\frac{(\xi - \eta^2)^{\frac{3}{2}}}{\eta^3} \right) \right] = -\frac{\xi}{\eta^3}$$

*Dato che $\psi_{\eta\eta}^{-1} = \psi_\xi^{-1} = 0$ vale $\psi_\xi^{-1} \phi_{\eta\eta}^{-1} = \phi_\xi^{-1} \psi_{\eta\eta}^{-1}$, per cui la scelta di $\psi(x,y)$ è compatibile con il requisito $\mathcal{L}^2[\psi] = 0$.*

*Es.2*

*Nell'esempio 2 del caso parabolico $\Phi$ e $\Phi^{-1}$ erano*

$$\begin{aligned} \phi(x,y) &= x + \frac{y^2}{2} = \xi & \phi^{-1}(\xi, \eta) &= \xi - \frac{\eta^2}{2} = x \\ \psi(x,y) &= y = \eta & \psi^{-1}(\xi, \eta) &= \eta = y \end{aligned}$$

*Bisogna fare però attenzione alla definizione di $\mathcal{L} = \partial_y - y\partial_x$ per cui : o si scambiano x ed y per poi invertirli nuovamente una volta trovato l'integrale generale, oppure si "ribalta" la trattazione. Procedendo per la seconda strada*

$$\mathcal{L}[\psi] = \psi_y - \Lambda\psi_x = \psi_y - \frac{\phi_y}{\phi_x}\psi_x = \frac{1}{\phi_x}(\phi_x\psi_y - \phi_y\psi_x) = \frac{|J_\Phi(x,y)|}{\phi_x} = \frac{1}{\psi_\eta^{-1}}$$

$$\mathcal{L} = \mathcal{L}[\psi]\partial_\eta = \frac{1}{\psi_\eta^{-1}}\partial_\eta$$

*Si ha quindi*

$$\mathcal{L}(\mathcal{L}[\psi]) = -\frac{\psi_{\eta\eta}^{-1}}{\left(\psi_\eta^{-1}\right)^3} \quad \Rightarrow \quad \mathcal{L}(\mathcal{L}[\psi]) = 0 \Leftrightarrow \psi_{\eta\eta}^{-1} = 0$$

*la quale è soddisfatta:* $\mathcal{L}[\psi] = (\partial_y - y\partial_x)y = 1 \Rightarrow \mathcal{L}(\mathcal{L}[\psi]) = \mathcal{L}[1] = 0$, *mentre le condizioni*

$$\mathcal{L}^2[\psi] = 0 \Leftrightarrow \psi_\xi^{-1}\phi_{\eta\eta}^{-1} = \phi_\xi^{-1}\psi_{\eta\eta}^{-1}$$

$$\mathcal{L}^2 = a\mathcal{L}\mathcal{L} \Leftrightarrow \mathcal{L}[\Lambda] = 0 \Leftrightarrow \phi_{\eta\eta}^{-1} = \psi_{\eta\eta}^{-1} = 0$$

*restano immutate. La prima, siccome la scelta di $\psi(x,y)$ è compatibile con il requisito $\mathcal{L}^2[\psi] = 0$, è soddisfatta, essendo $\psi_\xi^{-1} = \psi_{\eta\eta}^{-1} = 0$; mentre si ha che $\phi_{\eta\eta}^{-1} = -1 \neq 0$ e in effetti $\mathcal{L}[\Lambda] = 1 \neq 0$.*

*Si vuole ora proporre una nuova scelta (poco "felice" in quanto complica i calcoli) della funzione $\psi(x,y)$ al fine di mostrare che, nonostante la sua arbitrarietà (vincolata alla richiesta che $\mathcal{L}^2[\psi]$ sia nullo) il risultato non cambia. Sia dunque*

$$\phi(x,y) = x + \frac{y^2}{2} = \xi \qquad\qquad \phi^{-1}(\xi,\eta) = \eta - 1 - [1 + 2(\xi - \eta)]^{\frac{1}{2}} = x$$

$$\psi(x,y) = x + y = \eta \qquad\qquad \psi^{-1}(\xi,\eta) = 1 + [1 + 2(\xi - \eta)]^{\frac{1}{2}} = y$$

*Dato che $\mathcal{L}^2 = y^2\partial_{xx} - 2y\,\partial_{xy} + \partial_{yy}$, $\mathcal{L}^2[\psi] = 0$ in quanto $\psi_{xx} = \psi_{xy} = \psi_{yy} = 0$. Inoltre, dal momento che $\psi_y = 1$ ed $\mathcal{L}[\Lambda] = 1$, $\mathcal{L}(\mathcal{L}[\psi]) = -1 \neq 0$; in effetti da $\Phi^{-1}$,*

$$\phi_\xi^{-1} = -[1 + 2(\xi - \eta)]^{-\frac{1}{2}} \qquad \phi_\eta^{-1} = 1 + [1 + 2(\xi - \eta)]^{-\frac{1}{2}} \qquad \phi_{\eta\eta}^{-1} = [1 + 2(\xi - \eta)]^{-\frac{3}{2}}$$

$$\psi_\xi^{-1} = [1 + 2(\xi - \eta)]^{-\frac{1}{2}} \qquad \psi_\eta^{-1} = -[1 + 2(\xi - \eta)]^{-\frac{1}{2}} \qquad \psi_{\eta\eta}^{-1} = -[1 + 2(\xi - \eta)]^{-\frac{3}{2}}$$

*$\psi_{\eta\eta}^{-1} = -\phi_{\eta\eta}^{-1} \neq 0$ il che implica fra l'altro, in virtù dell'uguaglianza $\psi_\xi^{-1}\phi_{\eta\eta}^{-1} = \phi_\xi^{-1}\psi_{\eta\eta}^{-1}$ che anche $\phi_\xi^{-1} = -\psi_\xi^{-1} \neq 0$, come si può osservare. Tenendo sempre a mente che le formule viste in precedenza sono adattate al caso in cui*

$$\Lambda = \frac{\phi_y}{\phi_x} = -\frac{\phi_\eta^{-1}}{\psi_\eta^{-1}} = -\frac{1 + [1 + 2(\xi - \eta)]^{-\frac{1}{2}}}{-[1 + 2(\xi - \eta)]^{-\frac{1}{2}}} = [1 + 2(\xi - \eta)]^{\frac{1}{2}} + 1 = y$$

*si possono verificare*

$$\mathcal{L}(\mathcal{L}[\psi]) = -\frac{\psi_{\eta\eta}^{-1}}{\left(\psi_\eta^{-1}\right)^3} = -\frac{-[1+2(\xi-\eta)]^{-\frac{3}{2}}}{-[1+2(\xi-\eta)]^{-\frac{3}{2}}} = -1$$

$$\mathcal{L}[\Lambda] = \frac{\psi_{\eta\eta}^{-1}\phi_\eta^{-1} - \phi_{\eta\eta}^{-1}\psi_\eta^{-1}}{\left(\phi_\eta^{-1}\right)^3} = \frac{-[1+2(\xi-\eta)]^{-\frac{3}{2}} - [1+2(\xi-\eta)]^{-2} + [1+2(\xi-\eta)]^{-2}}{-[1+2(\xi-\eta)]^{-\frac{3}{2}}} = 1$$

*L'espressione per l'operatore parabolico stavolta sarà*

$$\mathcal{L}^2 = \Phi^{-1} \circ \{a \left[(\mathcal{L}[\psi])^2 \partial_{\eta\eta} + (\mathcal{L}[\Lambda]\phi_x)\partial_\xi\right]\}$$

*e considerato che $a = 1$, $\phi_x = 1$ e $(\mathcal{L}[\psi])^2 = \frac{1}{(\psi_\eta^{-1})^2} = 1 + 2(\xi - \eta)$ si deduce che*

$$\mathcal{L}^2 = [1 + 2(\xi - \eta)]\partial_{\eta\eta} + \partial_\xi$$

*L'equazione di partenza nelle variabili $(\xi, \eta)$ diventa*

$$y^2 u_{xx} - 2y\, u_{xy} + u_{yy} = u_x + 6y$$

$$[1 + 2(\xi - \eta)]U_{\eta\eta} + U_\xi = U_\xi + U_\eta + 6\{1 + [1 + 2(\xi - \eta)]^{\frac{1}{2}}\}$$

*ovvero un' ODE lineare del prim'ordine rispetto alla funzione $U_\eta$ nella variabile $\eta$*

$$[1 + 2(\xi - \eta)]U_{\eta\eta} - U_\eta = 6\{1 + [1 + 2(\xi - \eta)]^{\frac{1}{2}}\}$$

*Risolvendola si ottiene*

$$U_\eta = \frac{f(\xi) + 6\eta}{[1 + 2(\xi - \eta)]^{\frac{1}{2}}} - 6$$

*da cui integrando ancora una volta e sostituendo tenendo a mente che $f$ e $g$ sono funzioni arbitrarie*

$$U(\xi, \eta) = [1 + 2(\xi - \eta)]^{\frac{1}{2}} f(\xi) + g(\xi) - 2[1 + 2(\xi - \eta)]^{\frac{1}{2}}(2\xi + \eta + 1) - 6\eta$$

$$u(x, y) = y\, f\left(x + \frac{y^2}{2}\right) + g\left(x + \frac{y^2}{2}\right) - 2y^3 - 6xy - 6y + 2$$

*Sempre per la loro arbitrarietà, ponendo*

$$f\left(x + \frac{y^2}{2}\right) = F\left(x + \frac{y^2}{2}\right) + 6\left[\left(x + \frac{y^2}{2}\right) + 1\right] \qquad g\left(x + \frac{y^2}{2}\right) = G\left(x + \frac{y^2}{2}\right) - 2$$

*ci si riconduce a*

$$u(x, y) = y^3 + y\, F\left(x + \frac{y^2}{2}\right) + G\left(x + \frac{y^2}{2}\right)$$

*che è la forma che l'integrale generale della PDE assumeva nell'esempio 2.*

*Es.3*

*Poiché* $\mathcal{L}[\Lambda] = 0 \iff \phi_{\eta\eta}^{-1} = \psi_{\eta\eta}^{-1} = 0$ *ci si aspetta che in generale*

$$\begin{cases} \phi^{-1}(\xi,\eta) = \eta f(\xi) + g(\xi) \\ \psi^{-1}(\xi,\eta) = \eta p(\xi) + q(\xi) \end{cases}$$

*L'obiettivo è dare origine ad una generica equazione per cui* $\mathcal{L}^2 = a\mathcal{L}\mathcal{L}$ *. Si prenda ad esempio*

$$\begin{cases} \phi^{-1}(\xi,\eta) = c_1\eta\xi^n + c_2\xi^m - c_3 = x \\ \psi^{-1}(\xi,\eta) = c_4\eta\xi^{n-m} = y \end{cases}$$

*Raccogliendo* $\xi^m(c_1\eta\xi^{n-m} + c_2) = x + c_3$ *e sostituendo* $\eta\xi^{n-m} = \frac{y}{c_4}$ *si ottengono*

$$\xi = \left(\frac{c_4 x + c_3 c_4}{c_1 y + c_2 c_4}\right)^{\frac{1}{m}} \qquad \eta = \frac{y}{c_4}\left(\frac{c_4 x + c_3 c_4}{c_1 y + c_2 c_4}\right)^{1-\frac{n}{m}}$$

*Se* $[r(x,y)]^\alpha$ *è invariante lo sarà anche* $r(x,y)$ *per cui si può porre* $n = m = 1$

$$\phi(x,y) = \frac{c_4 x + c_3 c_4}{c_1 y + c_2 c_4} = \xi \qquad\qquad \phi^{-1}(\xi,\eta) = c_1\eta\xi + c_2\xi - c_3 = x$$
$$\psi(x,y) = \frac{y}{c_4} = \eta \qquad\qquad \psi^{-1}(\xi,\eta) = c_4\eta = y$$

*Si può ora ricavare facilmente* $\Lambda$

$$\Lambda = -\frac{\phi_x}{\phi_y} = -\frac{c_1 y + c_2 c_4}{c_1 x + c_1 c_3}$$

*Ribattezzando* $c_1 = c$ , $c_1 c_3 = a$ , $c_2 c_4 = b$ *si crea l'operatore*

$$\mathcal{L} = \left(\frac{\partial}{\partial x} + \frac{b + cy}{a + cx}\frac{\partial}{\partial y}\right)$$

*il cui prodotto* $\mathcal{L}^2 = (a + cx)^2 \mathcal{L}\mathcal{L}$ *applicato alla funzione incognita* $u(x,y)$, $\mathcal{L}^2[u] = 0$ *,genera l'equazione*

$$(a + cx)^2 u_{xx} + 2(a + cx)(b + cy)u_{xy} + (b + cy)^2 u_{yy} = 0$$

*con soluzione*

$$u(x,y) = (b + cy) F\left(\frac{a + cx}{b + cy}\right) + G\left(\frac{a + cx}{b + cy}\right)$$

*Per* $c = 0$ *si ottiene una generica PDE parabolica a coefficienti costanti*

$$a^2 u_{xx} + 2ab u_{xy} + b^2 u_{yy} = 0$$

*mentre per* $a = b = 0$ *si ritrova l'equazione*

$$x^2 u_{xx} + 2xy\, u_{xy} + y^2 u_{yy} = 0$$

*già vista nell'esempio 1.*

# *PDE II e trasformazioni del piano*

Un *diffeomorfismo* è una mappa tra varietà, differenziabile con inversa anch'essa differenziabile; questo termine è stato precedentemente usato come sinonimo di "cambiamento di carta" ovvero una trasformazione di coordinate. E' un concetto della geometria differenziale che richiama quello topologico di *omeomorfismo* (dal greco "forma-simile") e come tale può esprimere una "deformazione": si pensi alle curve caratteristiche $\phi(x,y) = k$, $\psi(x,y) = h$, le quali nel piano $(\xi, \eta)$ diventano le rette $\xi = k$, $\eta = h$, parallele agli assi $\eta$ e $\xi$ rispettivamente. Scopo di questo paragrafo è, partendo da $\Phi^{-1}$ e dalle trasformazioni (in generale non lineari) del piano, impostare un algoritmo per il passaggio in forma canonica di una PDE II, che si rivelerà poi equivalente al metodo "classico", con il fine di mettere in luce alcuni aspetti geometrici del problema. Si comincia con un esempio, per formalizzare i risultati nel seguito.

*Es.*

*Nell'esempio 5 del caso iperbolico si aveva*

$$u_{tt} + 4t u_{tx} + 3t^2 u_{xx} = 0$$

$$\phi(t,x) = x - \frac{t^2}{2} = \xi \qquad\qquad \phi^{-1}(\xi,\eta) = (\xi - \eta)^{\frac{1}{2}} = t$$

$$\psi(t,x) = x - \frac{3t^2}{2} = \eta \qquad\qquad \psi^{-1}(\xi,\eta) = \frac{3\xi - \eta}{2} = x$$

*Si traccino nel piano $(\xi, \eta)$ le curve $\phi^{-1}(\xi, \eta) = t = cost.$ e $\psi^{-1}(\xi, \eta) = x = cost.$, ovvero le rette $\eta = \xi - t^2$ (in blu) ed $\eta = 3\xi - 2x$ (in rosso) al variare dei parametri t e x: si hanno dei fasci di rette le cui intersezioni indicano la posizione dei punti $(t, x)$ nel piano $(\xi, \eta)$. Si ottiene quindi una rappresentazione "suggestiva" di come lo spazio venga deformato affinché le parabole $x = \frac{t^2}{2} + k$, $x = \frac{3t^2}{2} + h$ si trasformino nelle rette $\xi = k$, $\eta = h$.*

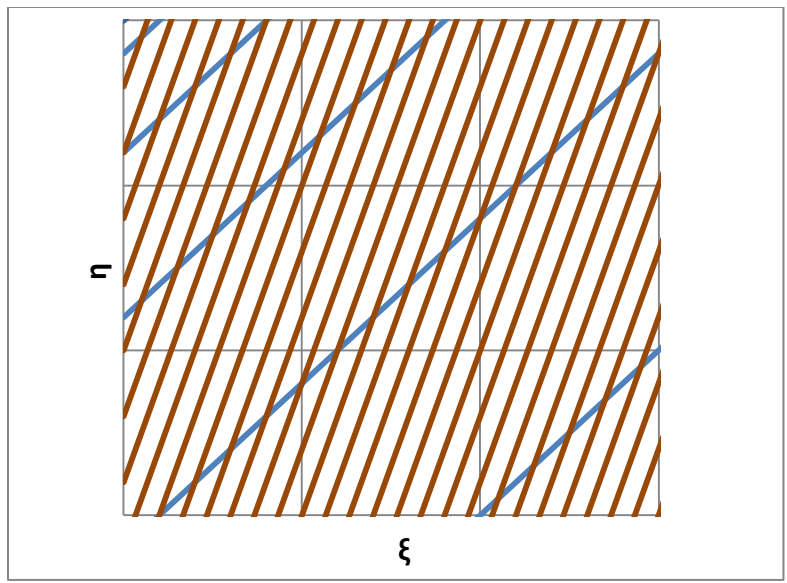

*Ci si proietti ora nel piano $(t,x)$: poiché derivare a $x = $ cost. equivale a derivare rispetto a $t$, parimenti derivare a $t = $ cost. significa derivare rispetto a $x$, nel piano $(\xi, \eta)$ derivare nella direzione di un vettore $\boldsymbol{v}$ tangente a $t = $ cost. (nell'esempio $\boldsymbol{v} = (1,1)$ ) sarà equivalente a derivare rispetto a $x$, mentre derivare nella direzione di un vettore $\boldsymbol{w}$ tangente a $x = $ cost. (qui $\boldsymbol{w} = (1,3)$ ) vorrà dire derivare rispetto a $t$. Ci si aspetta dunque le corrispondenze*

$$\frac{\partial}{\partial \boldsymbol{v}} \leftrightarrow \frac{\partial}{\partial x} \qquad \frac{\partial}{\partial \boldsymbol{w}} \leftrightarrow \frac{\partial}{\partial t}$$

*ricordando sempre che $\frac{\partial}{\partial \boldsymbol{w}} = \boldsymbol{w} \cdot \boldsymbol{\nabla} = \|\boldsymbol{w}\| \frac{\partial}{\partial \widehat{\boldsymbol{w}}}$ , $\|\widehat{\boldsymbol{w}}\| = 1$ . In effetti si ha proprio che*

$$\frac{\partial \phi^{-1}}{\partial \boldsymbol{v}} = 0 \qquad \frac{\partial \psi^{-1}}{\partial \boldsymbol{w}} = 0 \qquad \frac{\partial \psi^{-1}}{\partial \boldsymbol{v}} = 1 \qquad \frac{\partial \phi^{-1}}{\partial \boldsymbol{w}} = -(\xi - \eta)^{-\frac{1}{2}} = -\frac{1}{t}$$

*Adesso si applicano gli operatori $\frac{\partial}{\partial \boldsymbol{v}}, \frac{\partial}{\partial \boldsymbol{w}}$ alla funzione $u(t,x) = u(\phi^{-1}(\xi, \eta), \psi^{-1}(\xi, \eta))$*

$$\frac{\partial u}{\partial \boldsymbol{v}} = \frac{\partial u}{\partial t} \frac{\partial \phi^{-1}}{\partial \boldsymbol{v}} + \frac{\partial u}{\partial x} \frac{\partial \psi^{-1}}{\partial \boldsymbol{v}} = \frac{\partial u}{\partial x}$$

$$\frac{\partial u}{\partial \boldsymbol{w}} = \frac{\partial u}{\partial t} \frac{\partial \phi^{-1}}{\partial \boldsymbol{w}} + \frac{\partial u}{\partial x} \frac{\partial \psi^{-1}}{\partial \boldsymbol{w}} = -\frac{1}{t} \frac{\partial u}{\partial t}$$

*da cui*

$$\frac{\partial}{\partial t} = -(\xi - \eta)^{\frac{1}{2}} \frac{\partial}{\partial \boldsymbol{w}} = -(\xi - \eta)^{\frac{1}{2}} \left( \frac{\partial}{\partial \xi} + 3 \frac{\partial}{\partial \eta} \right)$$

$$\frac{\partial}{\partial x} = \frac{\partial}{\partial \boldsymbol{v}} = \frac{\partial}{\partial \xi} + \frac{\partial}{\partial \eta}$$

*Sapendo questo si calcolano*

$$\partial_{tt} = -(\xi - \eta)^{\frac{1}{2}}(\partial_\xi + 3\partial_\eta)\left[-(\xi - \eta)^{\frac{1}{2}}(\partial_\xi + 3\partial_\eta)\right] = (\xi - \eta)(\partial_{\xi\xi} + 6\partial_{\xi\eta} + 9\partial_{\eta\eta}) - \partial_\xi - 3\partial_\eta$$

$$\partial_{xx} = (\partial_\xi + \partial_\eta)(\partial_\xi + \partial_\eta) = \partial_{\xi\xi} + 2\partial_{\xi\eta} + \partial_{\eta\eta}$$

$$\partial_{tx} = -(\xi - \eta)^{\frac{1}{2}}(\partial_\xi + 3\partial_\eta)(\partial_\xi + \partial_\eta) = -(\xi - \eta)^{\frac{1}{2}}(\partial_{\xi\xi} + 4\partial_{\xi\eta} + 3\partial_{\eta\eta})$$

*A questo punto si può esprimere l'operatore $\mathcal{L} = \partial_{tt} + 4t\partial_{tx} + 3t^2\partial_{xx}$ nelle nuove coordinate*

$$\mathcal{L} = (\xi - \eta)(\partial_{\xi\xi} + 6\partial_{\xi\eta} + 9\partial_{\eta\eta}) - \partial_\xi - 3\partial_\eta - 4(\xi - \eta)(\partial_{\xi\xi} + 4\partial_{\xi\eta} + 3\partial_{\eta\eta}) + 3(\xi - \eta)(\partial_{\xi\xi} + 2\partial_{\xi\eta} + \partial_{\eta\eta})$$

$$\mathcal{L} = (\xi - \eta)\left[(1 - 4 + 3)\partial_{\xi\xi} + (6 - 16 + 6)\partial_{\xi\eta} + (9 - 12 + 3)\partial_{\eta\eta}\right] - \partial_\xi - 3\partial_\eta$$

$$\mathcal{L} = -4(\xi - \eta)\partial_{\xi\eta} - \partial_\xi - 3\partial_\eta$$

$$\mathcal{L}[u] = 0 \implies 4(\xi - \eta)U_{\xi\eta} + U_\xi + 3U_\eta = 0$$

*in accordo con quanto trovato più rapidamente mediante la nota espressione dell'operatore iperbolico nelle coordinate $(\xi, \eta)$. Dallo svolgimento si intuisce che, a partire da considerazioni di carattere geometrico, si è sostanzialmente arrivati a ridurre la PDE nella sua forma canonica applicando la "chain rule".*

Volendo generalizzare il procedimento, data una PDE del secondo ordine, a partire dall'equazione alle caratteristiche si trovano gli invarianti e quindi la trasformazione inversa di coordinate

$$\begin{cases} \phi^{-1}(\xi,\eta) = x \\ \psi^{-1}(\xi,\eta) = y \end{cases}$$

Si cercano i vettori tangenti alle curve definite implicitamente nel piano $(\xi,\eta)$ da $\phi^{-1}(\xi,\eta) = cost.$ ($x = cost.$) e $\psi^{-1}(\xi,\eta) = cost.$ ($y = cost.$), per cui, a imitazione degli operatori

$$\mathcal{L}^{\pm} = \boldsymbol{w}^{\pm} \cdot \boldsymbol{\nabla} = \frac{\partial}{\partial \boldsymbol{w}^{\pm}} \quad , \quad \boldsymbol{w}^{\pm} = (1, -\Lambda^{\pm}) \; , \quad \Lambda^{+} = \frac{\phi_x}{\phi_y} \; , \quad \Lambda^{-} = \frac{\psi_x}{\psi_y}$$

si possono prendere $\boldsymbol{v}$ e $\boldsymbol{w}$ tali che

$$\frac{\partial \phi^{-1}}{\partial \boldsymbol{v}} = \boldsymbol{\nabla}\phi^{-1} \cdot \boldsymbol{v} = (\phi_\xi^{-1}, \phi_\eta^{-1}) \cdot \left(1, -\frac{\phi_\xi^{-1}}{\phi_\eta^{-1}}\right) = 0$$

$$\frac{\partial \psi^{-1}}{\partial \boldsymbol{w}} = \boldsymbol{\nabla}\psi^{-1} \cdot \boldsymbol{w} = (\psi_\xi^{-1}, \psi_\eta^{-1}) \cdot \left(1, -\frac{\psi_\xi^{-1}}{\psi_\eta^{-1}}\right) = 0$$

e dunque

$$\frac{\partial}{\partial \boldsymbol{v}} = \partial_\xi - \frac{\phi_\xi^{-1}}{\phi_\eta^{-1}} \partial_\eta \qquad \frac{\partial}{\partial \boldsymbol{w}} = \partial_\xi - \frac{\psi_\xi^{-1}}{\psi_\eta^{-1}} \partial_\eta$$

Applicando gli operatori a $u(x,y) = u(\phi^{-1}(\xi,\eta), \psi^{-1}(\xi,\eta))$

$$\frac{\partial u}{\partial \boldsymbol{v}} = \frac{\partial u}{\partial x} \frac{\partial \phi^{-1}}{\partial \boldsymbol{v}} + \frac{\partial u}{\partial y} \frac{\partial \psi^{-1}}{\partial \boldsymbol{v}} = \frac{\partial \psi^{-1}}{\partial \boldsymbol{v}} \frac{\partial u}{\partial y}$$

$$\frac{\partial u}{\partial \boldsymbol{w}} = \frac{\partial u}{\partial x} \frac{\partial \phi^{-1}}{\partial \boldsymbol{w}} + \frac{\partial u}{\partial y} \frac{\partial \psi^{-1}}{\partial \boldsymbol{w}} = \frac{\partial \phi^{-1}}{\partial \boldsymbol{w}} \frac{\partial u}{\partial x}$$

si stabiliscono le corrispondenze

$$\frac{\partial}{\partial x} = \frac{1}{\phi_w^{-1}} \frac{\partial}{\partial \boldsymbol{w}} \qquad \frac{\partial}{\partial y} = \frac{1}{\psi_v^{-1}} \frac{\partial}{\partial \boldsymbol{v}}$$

Si osserva fra l'altro che, nonostante in generale $\frac{\partial}{\partial \boldsymbol{v}} \frac{\partial}{\partial \boldsymbol{w}} \neq \frac{\partial}{\partial \boldsymbol{w}} \frac{\partial}{\partial \boldsymbol{v}}$, per l'uguaglianza delle derivate miste è sempre vero che $\left[\frac{1}{\psi_v^{-1}} \frac{\partial}{\partial \boldsymbol{v}}, \frac{1}{\phi_w^{-1}} \frac{\partial}{\partial \boldsymbol{w}}\right] = 0$. A questo punto si calcolano

$$\frac{\partial \phi^{-1}}{\partial \boldsymbol{w}} = \left(\partial_\xi - \frac{\psi_\xi^{-1}}{\psi_\eta^{-1}} \partial_\eta\right)\phi^{-1} = \frac{1}{\psi_\eta^{-1}} \left(\phi_\xi^{-1} \psi_\eta^{-1} - \psi_\xi^{-1} \phi_\eta^{-1}\right) = \frac{|J_{\Phi^{-1}}(\xi,\eta)|}{\psi_\eta^{-1}} = \frac{1}{\phi_x}$$

$$\frac{\partial \psi^{-1}}{\partial \boldsymbol{v}} = \left(\partial_\xi - \frac{\phi_\xi^{-1}}{\phi_\eta^{-1}} \partial_\eta\right)\psi^{-1} = \frac{1}{\phi_\eta^{-1}} \left(\phi_\xi^{-1} \psi_\eta^{-1} - \psi_\xi^{-1} \phi_\eta^{-1}\right) = \frac{|J_{\Phi^{-1}}(\xi,\eta)|}{\phi_\eta^{-1}} = \frac{1}{\phi_y}$$

così

$$\frac{1}{\phi_w^{-1}} = \phi_x \qquad \frac{1}{\psi_v^{-1}} = \phi_y \qquad \frac{\partial}{\partial \boldsymbol{v}} = \frac{1}{\phi_y} \frac{\partial}{\partial y} \qquad \frac{\partial}{\partial \boldsymbol{w}} = \frac{1}{\phi_x} \frac{\partial}{\partial x}$$

ma siccome $-\frac{\psi_\xi^{-1}}{\psi_\eta^{-1}} = \frac{\psi_x}{\phi_x}$ e $-\frac{\phi_\xi^{-1}}{\phi_\eta^{-1}} = \frac{\psi_y}{\phi_y}$ si ritrovano

$$\frac{\partial}{\partial x} = \phi_x \left( \partial_\xi + \frac{\psi_x}{\phi_x} \partial_\eta \right) = \phi_x \partial_\xi + \psi_x \partial_\eta = \frac{\partial \Phi}{\partial x} \cdot \boldsymbol{\nabla}_{\xi,\eta}$$

$$\frac{\partial}{\partial y} = \phi_y \left( \partial_\xi + \frac{\psi_y}{\phi_y} \partial_\eta \right) = \phi_y \partial_\xi + \psi_y \partial_\eta = \frac{\partial \Phi}{\partial y} \cdot \boldsymbol{\nabla}_{\xi,\eta}$$

che sono le espressioni ottenute applicando gli operatori $\partial_x$ e $\partial_y$ a $U(\phi(x,y),\psi(x,y)) = U(\xi,\eta)$. Si può procedere al calcolo delle derivate di ordine superiore

$$\partial_{xx} = \Phi^{-1} \circ (\Phi_x \cdot \boldsymbol{\nabla}_{\xi,\eta})[\Phi^{-1} \circ (\Phi_x \cdot \boldsymbol{\nabla}_{\xi,\eta})] =$$

$$= \Phi^{-1} \circ [\Phi_x \cdot (\Phi_x \cdot \boldsymbol{\nabla}_{\xi,\eta})\boldsymbol{\nabla}_{\xi,\eta}] + \Phi^{-1} \circ (\Phi_x \cdot \boldsymbol{\nabla}_{\xi,\eta})[\Phi^{-1} \circ \Phi_x] \cdot \boldsymbol{\nabla}_{\xi,\eta}$$

Considerato che $\Phi_x \cdot (\Phi_x \cdot \boldsymbol{\nabla}_{\xi,\eta})\boldsymbol{\nabla}_{\xi,\eta} = (\Phi_x \cdot \boldsymbol{\nabla}_{\xi,\eta})^2$, $\Phi^{-1} \circ (\Phi_x \cdot \boldsymbol{\nabla}_{\xi,\eta}) = \{\partial_x\}_{(\xi,\eta)}$ (nelle coordinate $(\xi,\eta)$) e $\Phi^{-1} \circ \Phi_x = \{\Phi_x\}_{(\xi,\eta)}$ e dunque $\Phi^{-1} \circ (\Phi_x \cdot \boldsymbol{\nabla}_{\xi,\eta})[\Phi^{-1} \circ \Phi_x] = \{\partial_x \Phi_x\}_{(\xi,\eta)} = \Phi^{-1} \circ \Phi_{xx}$ si arriva a

$$\partial_{xx} = \Phi^{-1} \circ \left[ (\Phi_x \cdot \boldsymbol{\nabla}_{\xi,\eta})^2 + \Phi_{xx} \cdot \boldsymbol{\nabla}_{\xi,\eta} \right]$$

similmente

$$\partial_{yy} = \Phi^{-1} \circ \left[ (\Phi_y \cdot \boldsymbol{\nabla}_{\xi,\eta})^2 + \Phi_{yy} \cdot \boldsymbol{\nabla}_{\xi,\eta} \right]$$

Per il termine misto

$$\partial_{xy} = \Phi^{-1} \circ (\Phi_x \cdot \boldsymbol{\nabla}_{\xi,\eta})[\Phi^{-1} \circ (\Phi_y \cdot \boldsymbol{\nabla}_{\xi,\eta})] =$$

$$= \Phi^{-1} \circ [\Phi_y \cdot (\Phi_x \cdot \boldsymbol{\nabla}_{\xi,\eta})\boldsymbol{\nabla}_{\xi,\eta}] + \Phi^{-1} \circ (\Phi_x \cdot \boldsymbol{\nabla}_{\xi,\eta})[\Phi^{-1} \circ \Phi_y] \cdot \boldsymbol{\nabla}_{\xi,\eta}$$

ma $\Phi_y \cdot (\Phi_x \cdot \boldsymbol{\nabla}_{\xi,\eta})\boldsymbol{\nabla}_{\xi,\eta} = (\Phi_x \cdot \boldsymbol{\nabla}_{\xi,\eta})(\Phi_y \cdot \boldsymbol{\nabla}_{\xi,\eta})$ e $\Phi^{-1} \circ (\Phi_x \cdot \boldsymbol{\nabla}_{\xi,\eta})[\Phi^{-1} \circ \Phi_y] = \{\partial_x \Phi_y\}_{(\xi,\eta)} = \Phi^{-1} \circ \Phi_{xy}$, perciò

$$\partial_{xy} = \Phi^{-1} \circ \left[ (\Phi_x \cdot \boldsymbol{\nabla}_{\xi,\eta})(\Phi_y \cdot \boldsymbol{\nabla}_{\xi,\eta}) + \Phi_{xy} \cdot \boldsymbol{\nabla}_{\xi,\eta} \right]$$

Queste uguaglianze non sono nuove: sono già state menzionate precedentemente e si ottengono in forma "estesa" mediante la regola della catena, per poi compattarle "a mano". Per completezza, si vuole mostrare esplicitamente lo sviluppo al fine di mostrarne l'equivalenza; partendo da

$$u_x = \phi_x U_\xi + \psi_x U_\eta \qquad\qquad u_y = \phi_y U_\xi + \psi_y U_\eta$$

si calcola

$$u_{xx} = \partial_x(\phi_x U_\xi + \psi_x U_\eta)$$
$$= \phi_{xx} U_\xi + \psi_{xx} U_\eta + \phi_x\left[(U_\xi)_\xi \phi_x + (U_\xi)_\eta \psi_x\right] + \psi_x\left[(U_\eta)_\xi \phi_x + (U_\eta)_\eta \psi_x\right]$$
$$= \phi_x^2 U_{\xi\xi} + 2\phi_x \psi_x U_{\xi\eta} + \psi_x^2 U_{\eta\eta} + \phi_{xx} U_\xi + \psi_{xx} U_\eta$$

Per ottenere $u_{yy}$ è sufficiente scambiare gli indici nell'espressione precedente

$$u_{yy} = \phi_y^2 U_{\xi\xi} + 2\phi_y \psi_y U_{\xi\eta} + \psi_y^2 U_{\eta\eta} + \phi_{yy} U_\xi + \psi_{yy} U_\eta$$

mentre per il termine misto

$$u_{xy} = \partial_x\left(\phi_y U_\xi + \psi_y U_\eta\right)$$
$$= \phi_{xy} U_\xi + \psi_{xy} U_\eta + \phi_y\left[(U_\xi)_\xi \phi_x + (U_\xi)_\eta \psi_x\right] + \psi_y\left[(U_\eta)_\xi \phi_x + (U_\eta)_\eta \psi_x\right]$$
$$= \phi_x \phi_y U_{\xi\xi} + (\phi_x \psi_y + \phi_y \psi_x) U_{\xi\eta} + \psi_x \psi_y U_{\eta\eta} + \phi_{xy} U_\xi + \psi_{xy} U_\eta$$

Affinché i coefficienti di $U_{\xi\xi}$ e $U_{\eta\eta}$ siano nulli dovrà essere

$$a\phi_x^2 + 2b\phi_x\phi_y + c\phi_y^2 = 0 \qquad a\psi_x^2 + 2b\psi_x\psi_y + c\psi_y^2 = 0$$

Dividendo la prima per $\phi_y^2$ e la seconda per $\psi_y^2$, posti $\Lambda^+ = \frac{\phi_x}{\phi_y}$ e $\Lambda^- = \frac{\psi_x}{\psi_y}$, si ha che nel caso iperbolico

$$a(\Lambda^+)^2 + 2b\Lambda^+ + c = 0 \qquad a(\Lambda^-)^2 + 2b\Lambda^- + c = 0$$

oppure anche $a\mathcal{L}^-[\phi]\,\mathcal{L}^+[\phi] = 0$, $a\mathcal{L}^-[\psi]\,\mathcal{L}^+[\psi] = 0$, verificate in quanto $\mathcal{L}^+[\phi] = \mathcal{L}^-[\psi] = 0$.

Perché si annulli il coefficiente di $U_{\xi\eta}$

$$a\phi_x\psi_x + b(\phi_x\psi_y + \phi_y\psi_x) + c\phi_y\psi_y U_{\eta\eta}$$

Dividendo per $\phi_y\psi_y$ si ha

$$a\Lambda^+\Lambda^- + b(\Lambda^+ + \Lambda^-) + c = a\frac{c}{a} + b\left(-\frac{2b}{a}\right) + c = 2c - \frac{2b^2}{a} = 0 \;\Rightarrow\; b^2 = ac$$

cosa impossibile nel problema iperbolico dove $\Delta > 0$; invece ciò accade nel caso parabolico dove $\Delta = 0$, $a\phi_x^2 + 2b\phi_x\phi_y + c\phi_y^2 = a(\mathcal{L}[\phi])^2 = 0$ ma $a\psi_x^2 + 2b\psi_x\psi_y + c\psi_y^2 = a(\mathcal{L}[\psi])^2 \neq 0$.

Riassumendo, si è partiti da $\Phi^{-1} = (\phi^{-1}(\xi,\eta), \psi^{-1}(\xi,\eta))$ per trovare degli operatori $\partial_v$, $\partial_w$ nelle coordinate $(\xi,\eta)$ che avessero rispettivamente $\phi^{-1}$, $\psi^{-1}$ come invarianti, stabilendo delle relazioni tra $\partial_x$, $\partial_y$ e i suddetti operatori che permettono di esprimere direttamente $\partial_{xx}, \partial_{xy}, \partial_{yy}$ nelle nuove variabili. Applicando a $U(\xi,\eta)$ e sostituendo nell'equazione la si riduce in forma canonica. Il procedere *a ritroso* (nel senso che dagli invarianti si ricavano i relativi operatori e non viceversa) è una maniera per visualizzare la connessione tra il metodo delle caratteristiche e la geometria differenziale, in quanto dal punto di vista del mero calcolo è sostanzialmente equivalente all'usuale algoritmo risolutivo.

## *Conclusione*

Nel *problema iperbolico* l'operatore della parte principale $\mathcal{L}$ viene scomposto nel prodotto di due operatori $\mathcal{L}^+, \mathcal{L}^-$. Tuttavia, in generale, il prodotto $\mathcal{L}^- \mathcal{L}^+$ ( $\mathcal{L}^+ \mathcal{L}^-$ ) non è esattamente uguale a $\mathcal{L}$, ma comprende un termine aggiuntivo di ordine inferiore; tantomeno si può dire che il prodotto $\mathcal{L}^- \mathcal{L}^+$ sia sempre equivalente a $\mathcal{L}^+ \mathcal{L}^-$ (si hanno oramai diversi strumenti per riconoscere le varie casistiche). Comunque sia, trovando gli invarianti $\phi(x,y)$ e $\psi(x,y)$ a $\mathcal{L}^+, \mathcal{L}^-$ ( $\mathcal{L}^+[\phi] = 0$, $\mathcal{L}^-[\psi] = 0$ ) e cioè due funzioni le cui curve di livello $\phi(x,y) = k$, $\psi(x,y) = h$ abbiano come vettore tangente il vettore $\boldsymbol{v}^\pm$ associato al relativo operatore differenziale $\mathcal{L}^\pm = \boldsymbol{v}^\pm \cdot \boldsymbol{\nabla}$ (per cui una derivata direzionale sarebbe nulla) ed usandole come nuove coordinate $\phi(x,y) = \xi$, $\psi(x,y) = \eta$, si riesce a ridurre gli operatori differenziali $\mathcal{L}^\pm = \partial_x - \Lambda^\pm(x,y)\partial_y$ nelle variabili $(x,y)$ in semplici derivate parziali del prim'ordine nelle variabili $(\xi, \eta)$. Si è infatti visto che

$$\mathcal{L}^- = \frac{1}{\phi_\xi^{-1}} \partial_\xi \qquad \mathcal{L}^+ = \frac{1}{\phi_\eta^{-1}} \partial_\eta$$

dove la corrispondenza $\mathcal{L}^- \leftrightarrow \partial_\xi$ è giustificata dal fatto che si sta derivando lungo $\psi(x,y) = h$ ovvero $\eta = cost.$ nel piano $(\xi, \eta)$, mentre $\mathcal{L}^+ \leftrightarrow \partial_\eta$ in quanto si deriva lungo $\phi(x,y) = k$ e cioè $\xi = cost.$ nel piano $(\xi, \eta)$, ed effettivamente

$$\mathcal{L}^-[\psi] = \frac{1}{\phi_\xi^{-1}} \partial_\xi \eta = 0 \qquad \mathcal{L}^+[\phi] = \frac{1}{\phi_\eta^{-1}} \partial_\eta \xi = 0$$

E' chiaro quindi che, avendo la corrispondenza $\mathcal{L} \leftrightarrow \mathcal{L}^- \mathcal{L}^+$ ( $\mathcal{L}^+ \mathcal{L}^-$ ) (e non l'uguaglianza, essendo $\mathcal{L} = a \left( \mathcal{L}^- \mathcal{L}^+ + \mathcal{L}^-[\Lambda^+] \partial_y \right) = a \left( \mathcal{L}^+ \mathcal{L}^- + \mathcal{L}^+[\Lambda^-] \partial_y \right)$, dove in generale $\mathcal{L}^-[\Lambda^+]$ e $\mathcal{L}^+[\Lambda^-]$ non sono nulli), nelle coordinate $(\xi, \eta)$ si stabilirà la corrispondenza $\mathcal{L} \leftrightarrow \partial_{\xi\eta}$ (più eventuali termini di ordine inferiore dovuti a quanto detto poc'anzi), mentre non si presenteranno le derivate seconde $\partial_{\xi\xi}$ e $\partial_{\eta\eta}$, non avendo $\mathcal{L}^- \mathcal{L}^-$ o $\mathcal{L}^+ \mathcal{L}^+$.

Nel *problema parabolico* invece, $\mathcal{L}^2 = a(\mathcal{L}\mathcal{L} + \mathcal{L}[\Lambda]\partial_y)$ per cui si ha una sola famiglia di curve caratteristiche, soluzione di $\mathcal{L}[\phi] = 0$, quindi si completa il cambiamento di coordinate con una funzione "semiarbitraria" ( $\mathcal{L}[\psi] \neq 0$, $\mathcal{L}^2[\psi] = 0$ ). Siccome $\mathcal{L}[\phi] = 0$, $\phi(x,y) = cost.$ lungo la direzione definita dal vettore $\boldsymbol{v}$ in $\mathcal{L} = \boldsymbol{v} \cdot \boldsymbol{\nabla}$, dunque $\phi(x,y) = k$ saranno le sue linee di campo, quindi ci si aspetta la corrispondenza $\mathcal{L} \leftrightarrow \partial_\eta$ in quanto si deriva a $\xi = cost.$ nel piano $(\xi, \eta)$. E infatti si è trovato che

$$\mathcal{L} = \frac{1}{\phi_\eta^{-1}} \partial_\eta \implies \mathcal{L}[\phi] = \frac{1}{\phi_\eta^{-1}} \partial_\eta \xi = 0$$

ma allora se $\mathcal{L}^2 \leftrightarrow \mathcal{L}\mathcal{L}$, nelle variabili $(\xi, \eta)$ $\mathcal{L}^2 \leftrightarrow \partial_{\eta\eta}$ (più eventuali termini di ordine inferiore dovuti al fatto che in generale $\mathcal{L}[\Lambda] \neq 0$).

Le espressioni per l'operatore iperbolico e parabolico nelle coordinate $(\xi, \eta)$ scritte nella forma

$$\mathcal{L} = \Phi^{-1} \circ [a(\mathcal{L}^-[\phi]\,\mathcal{L}^+[\psi])\partial_{\xi\eta} + \mathcal{L}[\phi]\partial_\xi + \mathcal{L}[\psi]\partial_\eta] \qquad \mathcal{L}^2 = \Phi^{-1} \circ [a(\mathcal{L}[\psi])^2\, \partial_{\eta\eta} + \mathcal{L}^2[\phi]\partial_\xi]$$

riassumono questi concetti. Le uguaglianze presentano una certa compattezza, in quanto contengono al loro interno (a parte il coefficiente $a$ dell'equazione) soltanto gli operatori $\mathcal{L}^+, \mathcal{L}^-, \mathcal{L}$ e gli invarianti sulle caratteristiche $\phi(x,y)$ e $\psi(x,y)$. Danno già l'idea di quale potrà essere la forma canonica della parte principale, a seconda della natura del dominio (se iperbolico o parabolico) e del tipo di fattorizzazione dell'operatore principale, eludendo quindi il calcolo di coefficienti di derivate che si annulleranno. Offrono di conseguenza un algoritmo più leggero, più consapevole e meno macchinoso rispetto al calcolo di singoli termini (vedi le espressioni per $u_{xx}, u_{xy}, u_{yy}$ nel paragrafo precedente) che dovranno poi essere sostituiti all'interno della PDE.

Si termina con un esempio di riepilogo in cui si mettono a confronto, *step-by-step*, le diverse procedure per ridurre una PDE del secondo ordine in forma canonica ed eventualmente trovarne l'integrale generale.

*Es.*

$$xu_{xx} + (x-y)u_{xy} - yu_{yy} = 0$$

*Si comincia con il trovare*

$$x(\Lambda^\pm)^2 + (x-y)\Lambda^\pm - y = 0 \quad \Rightarrow \quad (\Lambda^\pm)^2 - \left(\frac{y}{x} - 1\right)\Lambda^\pm - \frac{y}{x} = 0 \quad \Rightarrow \quad \Lambda^+ = \frac{y}{x}, \quad \Lambda^- = -1$$

*Il problema è di tipo iperbolico. Da* $\frac{dy}{dx} = -\Lambda^\pm$

$$\frac{dy}{dx} = -\frac{y}{x} \quad \Rightarrow \quad \int \frac{dy}{y} = -\int \frac{dx}{x} \quad \Rightarrow \quad \ln y = -\ln x + cost. \quad \Rightarrow \quad \ln xy = cost. \quad \Rightarrow \quad xy = cost.$$

$$\frac{dy}{dx} = 1 \quad \Rightarrow \quad \int dy = \int dx \quad \Rightarrow \quad y = x + cost. \quad \Rightarrow \quad x - y = cost.$$

*si trovano gli invarianti*

$$\begin{cases} \phi(x,y) = xy &= \xi \\ \psi(x,y) = x - y &= \eta \end{cases}$$

*Questa procedura è comune ai vari metodi, in quanto si tratta di trovare le due famiglie di caratteristiche. Da qui si possono seguire diverse strade.*

*METODO I-A*

*Conoscendo* $\Lambda^\pm$ *si hanno gli operatori* $\mathcal{L}^\pm = \partial_x - \Lambda^\pm \partial_y$

$$\mathcal{L}^+ = \partial_x - \frac{y}{x}\partial_y \qquad \mathcal{L}^- = \partial_x + \partial_y$$

*Essendo* $\Lambda^-$ *un numero,* $\mathcal{L}^+[\Lambda^-] = 0$ *, mentre*

$$\mathcal{L}^-[\Lambda^+] = (\partial_x + \partial_y)\frac{y}{x} = -\frac{y}{x^2} + \frac{1}{x} = \frac{x-y}{x^2}$$

*Si osserva che* $\mathcal{L}^+[\Lambda^-] \neq \mathcal{L}^-[\Lambda^+]$ *e quindi gli operatori non commutano. E' però immediato sapere che* $\phi_y = x$, $\psi_y = -1$ *; si calcolano ancora*

$$\mathcal{L}^-[\phi] = (\partial_x + \partial_y)\,xy = y + x \qquad \mathcal{L}^+[\psi] = \left(\partial_x - \frac{y}{x}\partial_y\right)(x-y) = 1 + \frac{y}{x} = \frac{x+y}{x}$$

$$\mathcal{L}^-[\phi]\mathcal{L}^+[\psi] = \frac{(x+y)^2}{x}$$

*Non resta che inserire il tutto nell'espressione*

$$\mathcal{L} = \Phi^{-1}\circ\{a\,[\,(\mathcal{L}^-[\phi]\,\mathcal{L}^+[\psi])\,\partial_{\xi\eta} + (\mathcal{L}^-[\Lambda^+]\,\phi_y)\,\partial_\xi + (\mathcal{L}^+[\Lambda^-]\,\psi_y)\,\partial_\eta\,]\}$$

$$\mathcal{L} = \Phi^{-1}\circ\left\{x\left[\frac{(x+y)^2}{x}\partial_{\xi\eta} + \frac{x-y}{x^2}x\,\partial_\xi + 0\cdot(-1)\,\partial_\eta\right]\right\}$$

$$\mathcal{L} = \Phi^{-1}\circ\{[(x-y)^2 + 4xy]\,\partial_{\xi\eta} + (x-y)\,\partial_\xi\}$$

$$\mathcal{L} = (\eta^2 + 4\xi)\partial_{\xi\eta} + \eta\partial_\xi$$

*e applicando $\mathcal{L}$ a $U(\xi,\eta)$ si trova la forma canonica*

$$(\eta^2 + 4\xi)U_{\xi\eta} + \eta U_\xi = 0$$

*METODO I-B*

*Come il precedente, usando la seconda formulazione dell'espressione per l'operatore iperbolico nelle nuove coordinate; allora servono*

$$\phi_x = y \qquad \phi_{xx} = 0 \qquad \phi_y = x \qquad \phi_{yy} = 0 \qquad \phi_{xy} = 1$$

$$\mathcal{L}[\phi] = x\phi_{xx} + (x-y)\phi_{xy} - y\phi_{yy} = x - y$$

*$\psi(x,y) = x - y$ è lineare nelle sue variabili e quindi le derivate seconde si annullano; di conseguenza*

$$\mathcal{L}[\psi] = x\psi_{xx} + (x-y)\psi_{xy} - y\psi_{yy} = 0$$

*Ciò significa fra l'altro che qualsiasi funzione di $\psi$ è soluzione particolare dell'equazione differenziale. Si calcolano come nel precedente $\mathcal{L}^-[\phi]$ e $\mathcal{L}^+[\psi]$ e il prodotto $\mathcal{L}^-[\phi]\mathcal{L}^+[\psi]$ che è comodo scrivere come*

$$\mathcal{L}^-[\phi]\mathcal{L}^+[\psi] = \frac{(x-y)^2 + 4xy}{x}$$

*Si inserisce tutto in*

$$\mathcal{L} = \Phi^{-1}\circ[a(\mathcal{L}^-[\phi]\,\mathcal{L}^+[\psi])\partial_{\xi\eta} + \mathcal{L}[\phi]\partial_\xi + \mathcal{L}[\psi]\partial_\eta]$$

$$\mathcal{L} = \Phi^{-1}\circ\left[x\frac{(x-y)^2 + 4xy}{x} + (x-y)\partial_\xi + 0\cdot\partial_\eta\right]$$

*e come prima si ottiene*

$$\mathcal{L} = (\eta^2 + 4\xi)\partial_{\xi\eta} + \eta\partial_\xi$$

*METODO II*

*Da*

$$\begin{cases} \phi(x,y) = xy & = \xi \\ \psi(x,y) = x - y & = \eta \end{cases}$$

*invertendo si trova* $\Phi^{-1} = (\phi^{-1}(\xi,\eta), \psi^{-1}(\xi,\eta))$

$$x = \frac{\xi}{y} \;\Rightarrow\; \frac{\xi}{y} - y = \eta \;\Rightarrow\; y^2 + \eta y - \xi = 0 \;\Rightarrow\; y = \frac{1}{2}\left(\sqrt{\eta^2 + 4\xi} - \eta\right) = \psi^{-1}(\xi,\eta)$$

$$x = \frac{\xi}{y} = \frac{2\xi}{\sqrt{\eta^2+4\xi}-\eta} = \frac{2\xi}{\sqrt{\eta^2+4\xi}-\eta} \cdot \frac{\sqrt{\eta^2+4\xi}+\eta}{\sqrt{\eta^2+4\xi}+\eta} = \frac{1}{2}\left(\sqrt{\eta^2+4\xi}+\eta\right) = \phi^{-1}(\xi,\eta)$$

*Dato che* $\phi^{-1}(\xi,\eta) = f(\xi,\eta) + \frac{\eta}{2}$ *e* $\psi^{-1}(\xi,\eta) = f(\xi,\eta) - \frac{\eta}{2}$, $\phi_\xi^{-1} = f_\xi = \psi_\xi^{-1}$ *e quindi anche senza un calcolo esplicito si capisce che* $\phi_{\xi\eta}^{-1} = \psi_{\xi\eta}^{-1}$ *e si osserva che*

$$\frac{\psi_{\xi\eta}^{-1}}{\phi_{\xi\eta}^{-1}} = 1 = -(-1) = -\Lambda^- \;\Rightarrow\; \mathcal{L}^+[\Lambda^-] = 0 \;\Rightarrow\; \mathcal{L} = x\mathcal{L}^+\mathcal{L}^-$$

*Si calcola direttamente*

$$\mathcal{L} = x\mathcal{L}^+\mathcal{L}^- = \frac{\phi^{-1}}{\phi_\xi^{-1}\phi_\eta^{-1}}\left(\partial_{\xi\eta} - \phi_{\xi\eta}^{-1}\mathcal{L}^-\right)$$

*dove*

$$\phi_\xi^{-1} = \frac{1}{2}\frac{4}{2\sqrt{\eta^2+4\xi}} = \frac{1}{\sqrt{\eta^2+4\xi}} \;\Rightarrow\; \mathcal{L}^- = \frac{1}{\phi_\xi^{-1}}\partial_\xi = \sqrt{\eta^2+4\xi}\,\partial_\xi$$

$$\phi_\eta^{-1} = \frac{1}{2}\left(\frac{2\eta}{2\sqrt{\eta^2+4\xi}} + 1\right) = \frac{\sqrt{\eta^2+4\xi}+\eta}{2\sqrt{\eta^2+4\xi}}$$

$$\phi_{\xi\eta}^{-1} = -\frac{1}{2}(\eta^2+4\xi)^{-\frac{3}{2}}(2\eta) = -\frac{\eta}{(\eta^2+4\xi)\sqrt{\eta^2+4\xi}} \;\Rightarrow\; -\phi_{\xi\eta}^{-1}\mathcal{L}^- = \frac{\eta}{(\eta^2+4\xi)}\partial_\xi$$

$$\frac{\phi^{-1}}{\phi_\xi^{-1}\phi_\eta^{-1}} = \frac{1}{2}\left(\sqrt{\eta^2+4\xi}+\eta\right)\sqrt{\eta^2+4\xi}\,\frac{2\sqrt{\eta^2+4\xi}}{\sqrt{\eta^2+4\xi}+\eta} = \eta^2+4\xi$$

*Si conclude sostituendo nell'espressione*

$$\mathcal{L} = (\eta^2+4\xi)\left(\partial_{\xi\eta} + \frac{\eta}{(\eta^2+4\xi)}\partial_\xi\right) = (\eta^2+4\xi)\partial_{\xi\eta} + \eta\partial_\xi$$

*METODO III*

*Si utilizzano le espressioni ricavate mediante la chain rule (o le trasformazioni del piano)*

$$\partial_{xx} = \Phi^{-1} \circ \left[(\Phi_x \cdot \boldsymbol{\nabla}_{\xi,\eta})^2 + \Phi_{xx} \cdot \boldsymbol{\nabla}_{\xi,\eta}\right]$$

$$\partial_{xy} = \Phi^{-1} \circ \left[(\Phi_x \cdot \boldsymbol{\nabla}_{\xi,\eta})(\Phi_y \cdot \boldsymbol{\nabla}_{\xi,\eta}) + \Phi_{xy} \cdot \boldsymbol{\nabla}_{\xi,\eta}\right]$$

$$\partial_{yy} = \Phi^{-1} \circ \left[(\Phi_y \cdot \boldsymbol{\nabla}_{\xi,\eta})^2 + \Phi_{yy} \cdot \boldsymbol{\nabla}_{\xi,\eta}\right]$$

*Servono quindi, a partire da*

$$\begin{cases} \phi(x,y) = xy & = \xi \\ \psi(x,y) = x - y & = \eta \end{cases}$$

$$\phi_x = y \quad \phi_{xx} = 0 \quad \phi_y = x \quad \phi_{yy} = 0 \quad \phi_{xy} = 1$$

$$\psi_x = 1 \quad \psi_{xx} = 0 \quad \psi_y = -1 \quad \psi_{yy} = 0 \quad \psi_{xy} = 0$$

*Per cui*

$$\Phi_x = (y,1) \quad \Phi_y = (x,-1) \quad \Phi_{xx} = \Phi_{yy} = (0,0) \quad \Phi_{xy} = (1,0)$$

*Andando a sostituire si ha*

$$\partial_{xx} = \Phi^{-1} \circ \left[\left((y,1) \cdot (\partial_\xi, \partial_\eta)\right)^2 + (0,0) \cdot (\partial_\xi, \partial_\eta)\right] = \Phi^{-1} \circ \left[(y\partial_\xi + \partial_\eta)^2\right]$$
$$= \Phi^{-1} \circ [y^2 \partial_{\xi\xi} + 2y\partial_{\xi\eta} + \partial_{\eta\eta}]$$

$$\partial_{xy} = \Phi^{-1} \circ \left[\left((y,1) \cdot (\partial_\xi, \partial_\eta)\right)\left((x,-1) \cdot (\partial_\xi, \partial_\eta)\right) + (1,0) \cdot (\partial_\xi, \partial_\eta)\right]$$
$$= \Phi^{-1} \circ [(y\partial_\xi + \partial_\eta)(x\partial_\xi - \partial_\eta) + \partial_\xi] = \Phi^{-1} \circ [xy\partial_{\xi\xi} + (x-y)\partial_{\xi\eta} - \partial_{\eta\eta} + \partial_\xi]$$

$$\partial_{yy} = \Phi^{-1} \circ \left[\left((x,-1) \cdot (\partial_\xi, \partial_\eta)\right)^2 + (0,0) \cdot (\partial_\xi, \partial_\eta)\right] = \Phi^{-1} \circ \left[(x\partial_\xi - \partial_\eta)^2\right]$$
$$= \Phi^{-1} \circ [x^2 \partial_{\xi\xi} - 2x\partial_{\xi\eta} + \partial_{\eta\eta}]$$

*Tenendo presente* $\mathcal{L} = x\partial_{xx} + (x-y)\partial_{xy} - y\partial_{yy}$ *si calcolano*

$$x\partial_{xx} = \Phi^{-1} \circ [xy^2 \partial_{\xi\xi} + 2xy\partial_{\xi\eta} + x\partial_{\eta\eta}]$$

$$(x-y)\partial_{xy} = \Phi^{-1} \circ [(x^2y - xy^2)\partial_{\xi\xi} + (x-y)^2 \partial_{\xi\eta} + (y-x)\partial_{\eta\eta} + (x-y)\partial_\xi]$$

$$-y\partial_{yy} = \Phi^{-1} \circ [-x^2 y \partial_{\xi\xi} + 2xy\partial_{\xi\eta} - y\partial_{\eta\eta}]$$

*Infine sommando si arriva a*

$$\mathcal{L} = \Phi^{-1} \circ [(xy^2 + x^2 y - xy^2 - x^2 y)\partial_{\xi\xi} + (2xy + (x-y)^2 + 2xy)\partial_{\xi\eta} + (x + y - x - y)\partial_{\eta\eta} + (x-y)\partial_\xi]$$

$$\mathcal{L} = \Phi^{-1} \circ [((x-y)^2 + 4xy)\partial_{\xi\eta} + (x-y)\partial_\xi] = (\eta^2 + 4\xi)\partial_{\xi\eta} + \eta\partial_\xi$$

*INEGRALE GENERALE*

*Si è visto che la forma canonica della PDE è*

$$(\eta^2 + 4\xi)U_{\xi\eta} + \eta U_\xi = 0$$

*Da qui, per completezza, se ne vuole trovare la soluzione. Allora*

$$(\eta^2 + 4\xi)\frac{dU_\xi}{d\eta} = -\eta U_\xi \quad \Rightarrow \quad \int \frac{dU_\xi}{U_\xi} = -\frac{1}{2}\int \frac{2\eta}{\eta^2 + 4\xi}\, d\eta \quad \Rightarrow \quad \ln U_\xi = -\frac{1}{2}\ln(\eta^2 + 4\xi) + f(\xi)$$

$$\ln[U_\xi (\eta^2 + 4\xi)^{\frac{1}{2}}] = f(\xi) \quad \Rightarrow \quad U_\xi(\eta^2+4\xi)^{\frac{1}{2}} = F(\xi) \quad \Rightarrow \quad \frac{dU}{d\xi} = \frac{F(\xi)}{(\eta^2+4\xi)^{\frac{1}{2}}} \quad \Rightarrow \quad \int dU = \int \frac{F(\xi)d\xi}{(\eta^2+4\xi)^{\frac{1}{2}}}$$

*per cui*

$$U(\xi,\eta) = \int \frac{F(\xi)}{\sqrt{\eta^2 + 4\xi}}\, d\xi\ +\ G(\eta)$$

*con $\xi = xy$, $\eta = x - y$.*

*Si osserva che in questo caso non è possibile scrivere l'integrale generale nelle variabili $(x,y)$; tuttavia si possono ottenere delle "famiglie" di soluzioni particolari fissando $F(\xi)$. Posto, ad esempio,*

$$F(\xi) = 2k \qquad k \in \mathbb{R}$$

*si risolve l'integrale*

$$2k\int (\eta^2 + 4\xi)^{-\frac{1}{2}}\, d\xi = \frac{k}{2}\int (\eta^2+4\xi)^{-\frac{1}{2}}\, d(4\xi) = k(\eta^2 + 4\xi)^{\frac{1}{2}}$$

*da cui sostituendo*

$$k(\eta^2 + 4\xi)^{\frac{1}{2}} = k[(x-y)^2 + 4xy]^{\frac{1}{2}} = k[(x+y)^2]^{\frac{1}{2}} = k(x+y)$$

*si trova*

$$u(x,y) = k(x+y) + G(x-y)$$